\documentclass[review,hidelinks,onefignum,onetabnum]{siamart220329}
\usepackage{lipsum,mathtools}
\usepackage{amsfonts}
\usepackage{graphicx}
\usepackage{epstopdf}
\usepackage{algorithmic}
\ifpdf
  \DeclareGraphicsExtensions{.eps,.pdf,.png,.jpg}
\else
  \DeclareGraphicsExtensions{.eps}
\fi

\newsiamremark{remark}{Remark}
\newsiamremark{hypothesis}{Hypothesis}
\crefname{hypothesis}{Hypothesis}{Hypotheses}
\newsiamthm{claim}{Claim}
\headers{CEM-GMsFEM for Schr\"{o}dinger equations}{X. Jin, L. Liu, X. Zhong and E. T. Chung}
\title{Efficient numerical method for the Schr{\"o}dinger equation with high-contrast potentials\thanks{Submitted to the editors February 8, 2025.
\funding{L. Liu acknowledges the support by National Key R\&D Program of China (2021YFA1001200), Ministry of Science and Technology in China, Early Career Scheme (24301021) and General Research Fund (14303022 \& 14301423) funded by Research Grants Council of Hong Kong. E. Chung's research is partially supported by the Hong Kong RGC General Research Fund Projects 14305423 and 14305624.}}}
\author{Xingguang Jin\thanks{Department of Mathematics, The Chinese University of Hong Kong, Hong Kong
(\email{xgjin@math.cuhk.edu.hk}).}
\and Liu Liu\thanks{Department of Mathematics, The Chinese University of Hong Kong, Hong Kong
(\email{lliu@math.cuhk.edu.hk}).}
\and Xiang Zhong\thanks{Corresponding author. Department of Mathematics, The Chinese University of Hong Kong, Hong Kong
(\email{xzhong@math.cuhk.edu.hk}).}
\and Eric T. Chung\thanks{Department of Mathematics, The Chinese University of Hong Kong, Hong Kong
  (\email{tschung@math.cuhk.edu.hk}.})}

\usepackage{amsopn}

\ifpdf
\hypersetup{
  pdftitle={An Example Article},
  pdfauthor={D. Doe, P. T. Frank, and J. E. Smith}
}
\fi
\usepackage{amsmath}
\nolinenumbers
\usepackage[utf8]{inputenc}
\DeclareUnicodeCharacter{0308}{\"{}} 
\usepackage{bm} % Bold math font.
\usepackage{physics}
\usepackage{graphicx}
\usepackage{subcaption}
\newtheorem{assumption}{Assumption}
\usepackage{tikz}
\newcommand{\circled}[1]{\tikz[baseline=(char.base)]{
\node[shape=circle, draw, inner sep=0.5mm] (char) {#1};}}
\newcommand{\dx}{\,\mathrm{d}}

\begin{document}
\maketitle
% REQUIRED
\begin{abstract}
In this paper, we study the Schr\"{o}dinger equation in the semiclassical regime and with multiscale potential function. 
We develop the so-called constraint energy minimization generalized multiscale finite element method (CEM-GMsFEM), in the framework of Crank-Nicolson (CN) discretization in time. The localized multiscale basis functions are constructed by addressing the spectral problem and a constrained energy minimization problem related to the Hamiltonian norm. 
A first-order convergence in the energy norm and second-order convergence in the $L^2$ norm for our numerical scheme are shown, with a relation between oversampling number in the CEM-GMsFEM method, spatial mesh size and the semiclassical parameter provided. Furthermore, we demonstrate the convergence of the proposed Crank-Nicolson CEM-GMsFEM scheme. The convergence requires $H/\sqrt{\Lambda}=O(\varepsilon^{\frac{5}{4}})$, $\Delta t=O(\varepsilon^{\frac{5}{4}})$ if $\varepsilon\leq \delta$; while if $\delta<\varepsilon$, the convergence requires $H/\sqrt{\Lambda}=O(\varepsilon^{\frac{1}{4}}\delta)$, $\Delta t=O(\frac{\delta^2}{\varepsilon^{3/4}})$ (where $H$ represents the maximum diameter of coarse elements, $\Lambda$ is the minimal eigenvalue associated with the eigenvector not included in the auxiliary space, $\Delta t$ is the time step, $0 < \varepsilon\ll 1$ is the Planck constant and $\delta$ describes the multiscale structure of the potential).Several numerical examples including 1D and 2D in space, with high-contrast potential are conducted to demonstrate the efficiency and accuracy of our proposed scheme. 
\end{abstract}

% REQUIRED

% REQUIRED

\section{Introduction}
In quantum mechanics, the Schr\"{o}dinger equation serves as an important model to describe the behavior of quantum particles in materials with complex microstructures. In the semiclassical regime when $\varepsilon \ll 1$, it is well known that the Schr\"{o}dinger equation is in the high frequency regime, where the solution generates $O(\varepsilon)$ scaled oscillations in space and time. The high frequency of the solution, in addition to the multiscale potential we are considering in this paper, leads to unaﬀordable computational cost. Thus, developing efficient numerical method becomes a timely task. If one aims for direct simulation of the wave function, one of the best choices is the time splitting spectral method, by Bao, Jin,  Markowich and Sparber, see \cite{JMS-Acta, BJM2002}, where the meshing strategy $\Delta t=O(\varepsilon)$ and $\Delta x=O(\varepsilon)$ is required for moderate values of $\varepsilon$.
This meshing strategy is true if one wants to approximation the wave function. If only physical observables are to be approximated, then one can use $\Delta t=O(1). $   
%{In the regime of supercritical geometric optics, i.e., a severe restriction on the time step, one only needs $\Delta x=O(\varepsilon),\Delta t=o(\varepsilon).$} 

On the other hand, there are many approximation-type methods, which are valid in the limit $\varepsilon\to 0$, such as the level set method and the moment closure method based on the WKB analysis and the Wigner transform \cite{JMS-Acta}. 

Several wave packets-based methods have also been introduced, which reduce the full quantum dynamics to Gaussian wave packets dynamics \cite{hagedor199877,Heller76,Lee1982,altmann2020quantitative,henning2017crank}, and gain significant savings in the computational cost, such as the Gaussian beam method \cite{Jin_Wu_Yang_2011,qian2010} and 
the Frozen Gaussian beam method \cite{Kay1994,lu2012,lu2016,lu2018}. When the uncertainties are considered in the model, there are some works as well \cite{GWPT-UQ}. 

For these existing numerical methods, to the authors' best knowledge, Few studies have considered high-contrast potential functions with multi-scale external potentials. Recent advancements in nanotechnology have led to the fabrication of various material devices with customized functionalities. Examples include heterojunctions, which is often advantageous to engineer the electronic energy bands in many solid-state device applications, including semiconductor lasers \cite{gevaux2008}, solar cells \cite{louwen2016}, transistors and quantum metamaterials \cite{Quach11}. 
This can lead to high-contrast potentials, which are characterized by significant variations in potential energy across different regions. The finite element method (FEM) \cite{Dörfler1996} and finite difference method (FDM) \cite{Markowich1999} suffer from the heavily computational cost due the strong mesh size relations caused by the small Planck constant and multiscale structures of the potential. The well konwn time-splitting spectral method \cite{JMS-Acta} requires enough smoothness on both the potential and the initial condition to reach optimal accuracy; In the case of non-smooth potentials, especially high contrast potentials, the time-splitting spectral method would suffer
from reduced convergence order and great approximation errors.
% Multiscale methods.
The nature of heterogeneous potentials in Schr\"{o}dinger problems motivates the application of multiscale computational methods.  A multiscale method was first proposed in \cite{Chen2019} to efficiently solve the Schr\"{o}dinger equation with multiscale potential. They have since been applied to time-dependent \cite{chen2020} and random potentials \cite{chen2020multiscale}. The spirit of these methods was mainly motivated by the localized orthogonal decomposition (LOD) method \cite{Altmann2021, Maalqvist2021}, utilizing quasi-interpolation operators to decompose the solution into macroscopic and microscopic components \cite{Altmann2021,Maalqvist2021}. Recently, a novel multiscale method named Constraint Energy Minimization Generalized Multiscale Finite Element Method (CEM-GMsFEM) is initially developed by Chung, Efendiev, and Leung in \cite{Chung2018} which is aimed for the high-contrast problems and it has been successfully
applied to various partial differential equations arising from practical applications, see, e.g., \cite{jin2024,chung2023multiscale}. This approach and the LOD method share certain similarities. For instance, they both rely on the exponential decay and the oversampling technique introduced in \cite{Hou1997} and subsequently proved to improve convergence rates (ref.\ \cite{Efendiev2000,Efendiev2009,henning2014two}).
The Generalized Multiscale Finite Element Method (GMsFEM) was proposed by Efendiev, Galvis, and Hou in \cite{Efendiev2013}.
GMsFEMs leverage spectral decomposition to perform dimension reduction for the online space, exhibiting superior performance when dealing with high-contrast and channel-like coefficient profiles \cite{Chung2014a}.
``CEM'' is the acronym for ``Constraint Energy Minimizing''.
The novelty of CEM-GMsFEMs resides in replacing quasi-interpolation operators in LOD with element-wise eigenspace projections. We emphasize that the local spectral problems incorporate coefficient information. As a result, the exponential decay rates and local convergence are explicitly dependent on $\Lambda$, which is defined as the minimal eigenvalue corresponding to the eigenvector not included in the auxiliary space. Furthermore, $\Lambda$ remains stable across varying contrast ratios in numerical experiments.
Moreover, CEM-GMsFEMs introduce a relaxed version of the energy minimization problems to construct multiscale bases, which eliminates the necessity of solving saddle-point linear systems.
Our intention here is not to present a comprehensive review of multiscale computational methods from the community, and hence, notable advancements such as heterogeneous multiscale methods \cite{Weinan2003}, generalized finite element methods \cite{Babuska2011,Babuska2020,Ma2022}, and variational multiscale methods \cite{Hughes1995,Hughes2007,maalqvist2015computation,altmann2022localization,altmann2019localized,peterseim2024super,galvis2010domain} are not covered.

% Our contribution.
This article serves as an application of the CEM-GMsFEM to Schr\"{o}dinger problems. The CEM-GMsFEM, being a multiscale computational method, is specifically designed to handle the low regularity of the solution , and the construction of multiscale bases in this paper is tailored to the Schr\"{o}dinger setting with the Hamiltonian norm. 
For instance, the auxiliary space forms a core module in the original CEM-GMsFEM and is created by solving generalized eigenvalue problems.
It is worth noting that in building multiscale bases, we adhere to the relaxed version of the energy minimization problems in this paper, which offers implementation advantages.
We demonstrate that contrast robustness, which is an important feature of the original CEM-GMsFEM, is inherited in the proposed method. We first prove the global convergence under the CEM-GMsFEM framework. A first-order convergence in the energy norm and second-order convergence in the $L^2$ norm are proved in Theorem \ref{glothm} (see Remark \ref{re_glothm} for more details).
Then with several assumptions, we prove the existence of multiscale bases, the exponential decay property for the oversampling size, and the stability of the online space. By leveraging the global convergence property and the exponential decay of multiscale basis functions, along with Assumption \ref{oversampling layers}, we present an $O(\frac{H}{\sqrt{\Lambda}\varepsilon})$ convergence in the energy norm and an $O(\frac{H^2}{\Lambda\varepsilon^2})$ convergence in the $L^2$ norm for the spatial discretization. Next, we address the temporal regularity of our newly constructed Crank-Nicolson CEM-GMsFEM scheme. The final convergence theorem, presented in Theorem \ref{final convergence}, indicates that the error remains bounded by a small number when considering the relationship between $H$, $\Lambda$, $\Delta t$ and $\varepsilon,\delta$.
More precisely, the convergence requires $H/\sqrt{\Lambda}=O(\varepsilon^{\frac{5}{4}})$, $\Delta t=O(\varepsilon^{\frac{5}{4}})$ if $\varepsilon\leq \delta$; while if $\delta<\varepsilon$, the convergence requires $H/\sqrt{\Lambda}=O(\varepsilon^{\frac{1}{4}}\delta)$, $\Delta t=O(\frac{\delta^2}{\varepsilon^{3/4}})$ (see Remark \ref{relations between parameters} for more details). In this paper, we mainly address the challenges brought by the high contrast property of the potential, whereas consider moderately small Planck constant $\varepsilon$, but not too small. % Additionally, our error bound remains uniform with respect to $\varepsilon$.

%Organization
This paper is organized as follows. In Section \ref{sec2}, we introduce the model problem. The construction of multiscale bases in the proposed method is detailed in Section \ref{sec3}. All theoretical analysis for the proposed method is gathered in Section \ref{sec4} and Section \ref{sec5}.  To validate the performance of the proposed method, Section \ref{sec6} presents numerical experiments conducted on four different models. Finally, in Section \ref{sec7}, we conclude the paper.
\section{Preliminaries}\label{sec2}     
Throughout this paper, all functions are complex-valued and the conjugate of a function $v$ is denoted by $\overline{v}$. The spatial derivative is denoted by $D^{\bf{\mu}}$, where $D^{\bm{\mu}}v=\partial^{\mu_1}\cdot\cdot\cdot\partial^{\mu_d}$ with the multi-index $\bm{\mu}=(\mu_1,...,\mu_d)\in\mathbb{N}^d$ and $\abs{\bm{\mu}}=\mu_1+...+\mu_d$. The spatial $L^2$ inner product is denoted by $(\cdot,\cdot)$ with $(v,w)=\int_{\Omega}v\overline{w}$. The spatial $L^2$ norm is denoted by $\norm{\cdot}$ with $\norm{v}=(v,v)$. $\norm{\cdot}_{\infty}$ is the spatial $L^\infty$ norm with $\norm{v}_\infty=\text{ess sup}_{x\in\Omega}\abs{v(x)}$. We define $H^1_P(\Omega)=\{v\in V|v\text{ is periodic on }\partial\Omega\}$, where $\Omega$ is a bounded domain.  To simplify notation, we denote by $C$ a generic positive constant which may vary in different cases  but does not depend on the small parameters $\varepsilon,\delta$, the oversampling size $l$, the spatial mesh size $H$, and the time step size $\Delta t$. We consider the following Schr\"{o}dinger equation with multiscale potential in the semiclassical regime:
\begin{equation}
\label{model problem}
\begin{cases}
\mathrm{i}\varepsilon\partial_t u^{\varepsilon,\delta}=-\frac{1}{2}\varepsilon^2\Delta u^{\varepsilon,\delta}+V^\delta(x)u^{\varepsilon,\delta}\quad x\in\Omega,t\in(0,T]\\
u^{\varepsilon,\delta}\in H^1_P(\Omega)\quad t\in(0,T]\\
u^{\varepsilon,\delta}(0,x)=u^{\varepsilon}_0(x),\quad x\in\Omega
\end{cases}
\end{equation}
where $0 < \varepsilon\ll 1$ is the Planck constant, 
$u^{\varepsilon,\delta} = u^{\varepsilon,\delta}(t,x)$ is the electron wavefunction and $u^{\varepsilon}_0(x)$ is the initial data that is dependent on $\varepsilon$. %(motivated by the WKB approximation). 
Here $V^\delta(x)$ is an external and given multiscale potential depending on the small parameter $0\leq \delta\ll 1$, which describes the multiscale structure of the potential.    
\begin{assumption}
\label{assumption for u0 and potential}
Assume that $u^{\varepsilon}_0(x)$ satisfies $\norm{D^{\bm{\mu}}u^{\varepsilon}_0}\leq \frac{C}{\varepsilon^{\abs{\bm{\mu}}}}$. We assume the multiscale potential $V^\delta(x)$ is a high-contrast with $V_\text{min}\leq V^\delta(x)\leq V_\text{max}$ for all $x\in\Omega$. More precisely, ``high contrast" means $|V_\text{min}|\ll |V_\text{max}|$. Note that the multiscale potential $V^\delta(x)$ may depend on $\delta$--a parameter we use to characterize the contrast levels; however the specific dependence of $V^\delta(x)$ on $\delta$ is determined by specific applications. We also assume the regularity $\norm{D^{\bm{\mu}}V^\delta(x)}_\infty\leq \frac{C}{\delta^{\abs{\bm{\mu}}}}$, which is necessary in our analysis.
\end{assumption}
We introduce the bilinear form associated with the Schr\"{o}dinger operator given by 
$$\mathcal{H}=-\frac{1}{2}\varepsilon^2\Delta+V^\delta$$
as 
\begin{align}
\label{bilinear form a}
a(w,v)=\frac{1}{2}\varepsilon^2(\nabla w,\nabla v)+(V^\delta w,v).
\end{align}
We define
\[
\norm{v}_a=a(v,v)^{\frac{1}{2}}=(\frac{\varepsilon^2}{2}\norm{\nabla v}^2+(V^\delta v,v))^{\frac{1}{2}}.
\]
If the stationary problem with $\mathcal{H}$ as the differential operator
\begin{equation}
\label{stationary problem}
\begin{cases}
\mathcal{H}u=f,\quad x\in\Omega\\
u, D^{\bm{\sigma}}_{\bm{x}}u\text{ are periodic on } \partial\Omega,\quad \abs{\bm{\sigma}}=1, 
\end{cases}
\end{equation}
is considered, where periodic boundary conditions are prescribed and $f\in L^2(\Omega)$, the associated variational problem would be to find $u\in H^1_P(\Omega)$ such that
\begin{align}
\label{variational form for stationary equation}
a(u,v)=(f,v), \quad \forall v\in H^1_P(\Omega).
\end{align}
In terms of the Lax-Milgram theorem, the variational problem (\ref{variational form for stationary equation}) admits a unique solution  $u\in  H^1_P(\Omega)$.  

We will discuss the construction of multiscale basis functions in the next section. We consider $V_\text{ms}$ to be the space spanned by all multiscale basis functions. Then the multiscale solution $u_\text{ms}$ is denoted as the solution of the following problem: find $u_\text{ms}\in V_\text{ms}$ such that
\begin{equation}
\label{multiscale problem}
a(u_\text{ms},v)=(f,v) \quad \forall v\in V_\text{ms}.
\end{equation}
\section{Multiscale basis functions}\label{sec3}   
To simplify the notation, we define $V\coloneqq H^1(\Omega)$ and $V(S)\coloneqq H^1(S)$ for all $S\subset\Omega$. Let $V(K_j)$ be the snapshot space on each coarse grid
block $K_j$, and we use the method of the spectral problem to solve an eigenvalue problem on $K_j$: find eigenvalues $\lambda_j^{i}\in \mathbb{R}$ and basis functions $\phi_j^{i}\in V(K_j)$ such that for all $v\in V(K_j)$,
\begin{equation}\label{eig}
a_j(\phi_j^i,v)=\lambda_j^i s_j(\phi_j^i,v),\quad\forall v\in {V(K_j)},
\end{equation}
$$
\begin{aligned}
a_j(\phi_j^i,v)=&\int_{K_j}\frac{1}{2}\varepsilon^2\nabla\phi_j^i\cdot\nabla\overline v +V^\delta(x)\phi_j^i\cdot \overline{v}\dx x,\\
s_j(\phi_j^i,v)=&\int_{K_j}{\tilde{\mu}(x)}\phi_j^i\cdot\overline{v}\dx x,\quad \tilde{\mu}(x) :=\sum_{k=1}^{N_v}\frac{1}{2}\varepsilon^2\nabla\eta_{j,k}^{i}\cdot\nabla\eta_{j,k}^{i},
\end{aligned}
$$
where $N_v$ is the number of vertices contained in an element, to be specific, $N_v=4$ for a quadrilateral mesh and  $\*{\eta_{j,1}^{i},\eta_{j,2}^{i},\cdots,\eta_{j, N_v}^{i}}$ is the set of Lagrange basis functions on the coarse element $K_j\in\mathcal{T}_H.$ The bilinear form $s_j(\cdot,\cdot)$ above defines an inner product with norm $\norm{v}_{s(K_j)}=\sqrt{s_j(v,v)}$.

Let the eigenvalues $\lambda_j^i$ in the ascending order:
$$
0=\lambda_j^0<\lambda_j^1\leq\lambda_j^2\leq\cdots\lambda_j^{l_j+1}\leq\cdots ,
$$ 
and we use the first $l_j$ eigenvalue functions corresponding to the eigenvalues to construct the local auxiliary space  $V_{\text{aux}}^{j}=\{\phi_j^1,\phi_j^2,\cdots,\phi_j^{l_j}\}.$ 
% {In terms of the usual Poincar$\acute{\rm e}$ inequality and \cite[Section 3.3]{galvis2010domain}, we know the problem (\ref{eig}) has eigenvalues that scale as the inverse of the high contrast. The number of eigenvalues that scale as the inverse of high conductivity is the same as the number of connected high-conducting regions and the rest of the eigenvalues are large and remain bounded below as the contrast increases. Therefore, if $l_j$ is large, then $\lambda_j^{l_j+1}$ is contrast independent, which is useful for the following analysis.}
The global auxiliary space $V_{\text{aux}}$ is the sum of these local auxiliary spaces, namely 
$V_{\text{aux}}= \bigoplus_{j=1}^N V_{\text{aux}}^{j}$, which will be used to construct multiscale basis functions.
The next we give the definition of the so called $\phi_j^i$-orthogonal, for a given a function $\phi_j^i\in V_{\text{aux}}$, $\psi\in V$, and we define
$$
s(\phi_j^i,\psi)=1, \quad s(\phi_{j'}^{i'},\psi)=0 \,\,\text{if}\, j'\neq j\,\text{or}\,i'\neq i.
$$
Based on the $\phi_j^i$-orthogonal, we can obtain that for any $v\in V$
$$s(\phi_j^i,v)=\sum_{j=1}^{N}s_j(\phi_j^i,v).\quad $$
Then $\norm{v}_{s}=s(v,v)^{\frac{1}{2}}$ for all $v\in V$ is the induced norm derived from the bilinear form $s(\cdot,\cdot)$. The orthogonal projection $\pi_j$ from $V(K_j)$ onto $V_\text{aux}^{j}$ is
$$\pi_{j}(v) :=\sum_{i=1}^{l_j}\frac{s(\phi_j^i,v)}{s(\phi_j^i,\phi_j^i)}\phi_j^i, \quad \forall v\in V(K_j),$$
and the global projection is $\pi :=\sum_{j=1}^{N}\pi_j$ from $H^1$ to $V_\text{aux}$
\footnote{We use a zero-extension here, which extends each $V_{\text{aux}}^j$ into $L^2(\Omega)$. }
We can immediately derive the following \cref{inter}, which shows an important property of the global projection $\pi$. We define $s^{-1}$-norm as: $\norm{v}_{s^{-1}}=(\int_{\Omega}\tilde{\mu}^{-1}\abs{v}^2\dx x)^{1/2}$ for all $v\in V$. Clearly $\norm{\cdot}_{s^{-1}}=O(H/\varepsilon)$.

\begin{lemma}\label{inter}
In each $K_j\in\mathcal{T}_H$, for all $v\in V(K_j)$,
%c_*^{-1}H^{-2}\varepsilon^2\norm{v-\pi_j v}_{(K_j}^2\leq
\begin{equation}\label{pi1}
\norm{v-\pi_j v}_{s(K_j)}^2\leq \frac{\norm{v}_{a(K_j)}^2}{\lambda_{l_j+1}^j}\leq\Lambda^{-1}\norm{v}_{a(K_j)}^2.
\end{equation}
where $\Lambda=\min_{1\leq j\leq N}\lambda_j^{l_j+1}$, and  
\begin{equation}\label{pi2}
\norm{\pi_jv}_{s(K_j)}^2=\norm{v}_{s(K_j)}^2-\norm{v-\pi_jv}_{s(K_j)}^2\leq	\norm{v}_{s(K_j)}^2.
\end{equation}
\end{lemma}
For each coarse element $K_j\in\mathcal{T}_H$, the oversampling domain $K_j^m\subset\Omega$ is constructed by enlarging $K_j$ for $m$ coarse grid layers. We define the multiscale basis function: find $\psi_{j,m}^i\in V(K_j^m)$ such that
\begin{equation}\label{msbasis1}
a(\psi_{j,m}^i,v)+s(\pi\psi_{j,m}^i,\pi v)=s(\phi_j^i,\pi v),\quad\forall v\in V(K_j^m),
\end{equation}
where $V(K_j^m)$ is the restriction of $V$ in $K_j^m.$
Now, our multiscale finite element space $V_{\text{ms}}$ can be defined by solving an variational problem \cref{msbasis1}
$$
V_{\text{ms}}=\text{span}\{\psi_{j,m}^{i}\,|\, 1\leq i\leq l_j,\,1\leq j\leq N\}.
$$
The global multiscale basis function $\psi_j^i\in V$ is defined similarly, 
\begin{equation}\label{msbasis2}
a(\psi_{j}^i,v)+s(\pi\psi_{j}^i,\pi v)=s(\phi_{j}^i,\pi v),\quad\forall v\in V.
\end{equation}
Thereby, the global multiscale finite element space $V_{\text{glo}}$  is defined by
$$
V_{\text{glo}}=\text{span}\{\psi_j^i\,|\, 1\leq i\leq l_j,\,1\leq j\leq N\}.
$$
Then the global multiscale solution $u_\text{glo}$ is denoted as the solution of the following problem: find $u_\text{glo}\in V_\text{glo}$ such that
\begin{equation}
\label{global multiscale problem}
a(u_\text{glo},v)=(f,v) \quad \forall v\in V_\text{glo}.
\end{equation}
\begin{lemma}
For any $v\in V_{\text{glo}}$, then $a(v,v')=0$ for any $v'\in V$ with $\pi v'=0$. If there exists $v'\in V$ such that $a(v,v')=0$ for any $v\in V_{\text{glo}}$, then $\pi v'=0$.
\end{lemma}
\section{Analysis}\label{sec4}   
\begin{theorem}\label{glothm}
Let $u$ be the solution to the stationary problem (\ref{stationary problem}) and $u_{\text{glo}}$ be the solution to the global problem (\ref{global multiscale problem}). If $f\in L^2(\Omega)$, then 
\begin{equation}
\label{energy estimate for u-u_glo}
\norm{u-u_{\text{glo}}}_a\leq \frac{\norm{f}_{s^{-1}}}{\sqrt{\Lambda}}
\end{equation}
and
\begin{equation}
\label{s norm estimate for u-u_glo}
\norm{u-u_\text{glo}}_s\leq \frac{C}{\Lambda}\norm{f}_{s^{-1}},
\end{equation}
where $C>0$ is independent of $\Lambda,\varepsilon,\delta,H$.
\end{theorem}
\begin{remark}
\label{second order convergence for l2 norm}
Notice that by (\ref{s norm estimate for u-u_glo}), we actually have (let $\norm{\cdot}$ define the $L^2$ norm on $\Omega$)
\[
\norm{u-u_\text{glo}}\leq C\frac{H^2}{\Lambda\varepsilon^2}\norm{f}
\]
since $\norm{\cdot}_s=O(\varepsilon/H)$ and $\norm{\cdot}_{s^{-1}}=O(H/\varepsilon)$. That is, we can obtain the second order convergence for the global error estimate with respect to $L^2$ norm.
\end{remark}
\begin{proof} For the global error estimate in the energy norm, 
\begin{equation*}
\begin{aligned}
 \norm{u-u_{\text{glo}}}_a^2&=a(u-u_{\text{glo}},u-u_{\text{glo}})=a(u,u-u_{\text{glo}})\\
&=(f,u-u_{\text{glo}})\leq \norm{f}_{s^{-1}}\norm{u-u_{\text{glo}}}_s\leq \frac{\norm{f}_{s^{-1}}}{\sqrt{\Lambda}}\norm{u-u_{\text{glo}}}_a,
\end{aligned}
\end{equation*}
where the last inequality comes from the fact that $\pi(u-u_{\text{glo}})=0$ and the use of the \cref{inter}.

Let $w\in V$ be the solution of
\[
a(w,v)=(u-u_\text{glo},v)\text{ for all } v\in V,
\]
and 
\[
a(w_\text{glo},v_\text{g})=(u-u_\text{glo},v_\text{g})\text{ for all } v_\text{g}\in V_\text{glo}.
\]
Then we easily obtain
\begin{align*}
\norm{u-u_\text{glo}}_s^2&\leq C\frac{\varepsilon^2}{H^2}\norm{u-u_\text{glo}}^2=C\frac{\varepsilon^2}{H^2}a(w,u-u_\text{glo})=C\frac{\varepsilon^2}{H^2}a(w-u_\text{glo},u-u_\text{glo})\\
&\leq C\frac{\varepsilon^2}{H^2}\norm{w-u_\text{glo}}_a\cdot\norm{u-u_\text{glo}}_a\leq C\frac{\varepsilon^2}{\Lambda H^2}\norm{f}_{s^{-1}}\norm{u-u_\text{glo}}_{s^{-1}}\\
&\leq C\frac{\varepsilon^2}{\Lambda H^2}\norm{f}_{s^{-1}}\frac{H^2}{\varepsilon^2}\norm{u-u_\text{glo}}_s=C\frac{1}{\Lambda}\norm{f}_{s^{-1}}\norm{u-u_\text{glo}}_s.
\end{align*}
Thus we have
$
\norm{u-u_\text{glo}}_s\leq C/\Lambda\norm{f}_{s^{-1}}.
$
This completes the proof.
\end{proof}
\begin{remark}\label{re_glothm}
From Theorem \ref{glothm}, the upper bound of the global error does not explicitly depend on $\varepsilon$ since the term $H/\varepsilon$ is hidden in the norm $\norm{\cdot}_{s^{-1}}$. More precisely, since $\norm{\cdot}_{s^{-1}}=O(H/\varepsilon)$, by (\ref{energy estimate for u-u_glo}) and Remark \ref{second order convergence for l2 norm}, we have $\norm{u-u_{\text{glo}}}_a\leq \frac{CH\norm{f}}{\varepsilon\sqrt{\Lambda}}$, $\norm{u-u_\text{glo}}\leq C\frac{H^2}{\Lambda\varepsilon^2}\norm{f}$. In fact, we share the same convergence order as \cite{wu2022} (i.e. first-order convergence in the energy norm and second-order convergence in the $L^2$ norm). 
\end{remark}
By the above theorem, we obtain the convergence of the method for using global basis functions. Then we will prove these global basis functions are localizable. We first give a lemma in the following before estimating the difference between the global basis functions and the multiscale basis functions. For each coarse block $K$, we define $B$ to be a bubble function with with $B(x)>0$ for all $x$ belonging to the interior of $K$ and $B(x)=0$ for all $x\in\partial K$. We will take $B=\Pi_k\eta_k$, where the product is taken over all vertices $k$ on the boundary of $K$. Utilizing this bubble function, we define the constant $$C_\pi=\sup_{K\in\mathcal{T}_H, v\in V_\text{aux}}\frac{\int_K \tilde{\mu}(x)v^2}{\int_KB(x)\tilde{\mu}(x)v^2}.$$
\begin{lemma}
\label{C inverse}
For all $v_\text{aux}\in V_\text{aux}$, there exists a function $v\in V$ such that 
\[
\pi(v)=v_\text{aux},\quad \norm{v}_a^2\leq C_0\norm{v_\text{aux}}_s^2,\quad \text{supp}(v)\subset \text{supp}(v_\text{aux}),
\]
where $C_0=2(1+N_v)^2(1+\lambda_\text{max})^2$ and $\lambda_\text{max}=\max_{1\leq j\leq N}\max_{1\leq i\leq l_j}\lambda^i_j$.
\end{lemma}
\begin{proof}
Consider the following minimization problem defined on a coarse block $K_j$:
\begin{equation}
\label{minimization problem}
v=\text{argmin}\{a(\psi,\psi)|\,\psi\in V(K_j),\,s_j(\psi,v_\text{aux})=1,\,s_j(\psi,w)=0\,\text{for all} \, w\in v_\text{aux}^\perp.\}
\end{equation}
Here $v_\text{aux}^\perp\in V_\text{aux}^j$ is the orthogonal complement of $\text{span}\{v_\text{aux}\}$ with respect to the inner product $s_j$. Let $p\in V_\text{aux}^j$ The minimization problem (\ref{minimization problem}) is equivalent to the following variational problem: find $v\in V(K_j)$ and $y\in V_\text{aux}^j$ such that
\begin{subequations}
\label{equivalent variational problem}
\begin{align}
\int_{K_j}\frac{1}{2}\varepsilon^2\nabla v\cdot\nabla\overline{w}+V^\delta(x)v\cdot\overline{w}dx+\int_{K_j}\tilde{\mu}(x)w\cdot\overline{y}dx&=0\quad \forall w\in V(K_j) \label{equivalent variational problem_a}\\
\int_{K_j}\tilde{\mu}(x)v\cdot\overline{z}dx&=\int_{K_j}\tilde{\mu}(x)p\cdot\overline{z}dx\quad \forall z\in  V_\text{aux}^j.\label{equivalent variational problem_b}
\end{align}
\end{subequations}
Note that the well-posedness of the minimization problem (\ref{minimization problem}) is equivalent to the existence of a function $v\in V(K_j)$ such that 
\[
s_j(v,p)\geq C\norm{p}_{s(K_j)}^2,\quad \norm{v}_{a(K_j)}\leq C\norm{p}_{s(K_j)},
\]
where $C$ is a constant independent of the meshsize and the problem parameters.

Notice that $p$ is supported in $K_j$. We let $v=B(x)p$. Then it follows from the definition of $s_j$ that
\[
s_j(v,p)=\int_{K_j}\tilde{\mu}(x)B(x)p^2\geq C_\pi^{-1}\norm{p}_{s(K_j)}^2.
\]
Because $\nabla(B(x)p)=p\nabla B(x)+B(x)\nabla p$, $\abs{B}\leq1$ and $\abs{\nabla B(x)}^2\leq N_v^2\sum_k\abs{\nabla\eta_k}^2$, we have
$
\norm{B(x)p}_{a(K_j)}^2\leq 2(\norm{p}_{a(K_j)}^2+N_v^2\norm{p}_{s(K_j)}^2).
$
Then we can obtain
\[
\norm{v}_{a(K_j)}=\norm{B(x)p}_{a(K_j)}\leq \sqrt{2}\cdot(1+N_v)(\norm{p}_{a(K_j)}+\norm{p}_{s(K_j)}).
\]
By using the local spectral problem (\ref{eig}), we have
$
\norm{p}_{a(K_j)}\leq (\max_{1\leq i\leq l_j}\lambda^i_j)\norm{p}_{s(K_j)}.
$
Therefore, we finally have
\[
\norm{v}_{a(K_j)}\leq \sqrt{2}(1+N_v)(1+\max_{1\leq i\leq l_j}\lambda^i_j)\norm{p}_{s(K_j)},
\]
which proves the unique solvability of the minimization problem. $v,y$ satisfy (\ref{equivalent variational problem}). From (\ref{equivalent variational problem_b}), we can see that $\pi_j(v)=p$. Since $v=Bp$, it's clear that $\text{supp}(v)\subset\text{supp}(p)$. By the above estimates, we also have the desired estimate with $C_0=2(1+N_v)^2(1+\lambda_\text{max})^2$, where $\lambda_\text{max}=\max_{1\leq j\leq N}\max_{1\leq i\leq l_j}\lambda^i_j$. This completes the proof.
\end{proof}

To estimate the difference between the global and multiscale basis function, we need some notations for the oversampling domain and the cutoff function with respect to these oversampling domains. For each $K_j$, we recall that $K_j^m\subset\Omega$ is the oversampling coarse region by enlarging $K_j$ by $m$ oversampling size. For $M>m$, we define $\eta_j^{M,m}\in\text{span}\{\eta_k\}$ such that
\begin{subequations}
	\label{cutoff function}
	\begin{align}
\eta_j^{M,m}&\equiv1\quad {\rm in}\hspace{0.5em} K_j^m,\label{cutoff function_a}\\
\eta_j^{M,m}&\equiv0\quad {\rm in}\hspace{0.5em}\Omega\backslash K_j^M,\label{cutoff function_b}\\
0\leq\eta_j^{M,m}&\leq1\quad {\rm in}\hspace{0.5em}K_j^M\backslash K_j^m.\label{cutoff function_c}
\end{align}
\end{subequations}
Clearly we observe that $\frac{1}{2}\varepsilon^2\abs{\nabla\eta_j^{M,m}}^2\leq \tilde{\mu}(x)$.
Next we show that our multiscale basis functions have a decay property. We first give an assumption:
\begin{assumption}
\label{regularity for mesh}
There exists a positive constant $C_\text{ol}$ such that for all $K_j\in\mathcal{T}_H$ and $m>0$, 
\[
\#\{K\in \mathcal{T}_H|\, K\subset K_j^m\}\leq C_\text{ol}m^d.
\]
\end{assumption}

\begin{theorem}
\label{decay property}
Consider the oversampling domain $K_j^l$ with $l\geq1$. Let $\phi^i_j\in V_\text{aux}$ be a given auxiliary multiscale basis function. Let $\psi^i_{j,l}$ be the multiscale basis functions obtained in (\ref{msbasis1}) and let $\psi^i_j$ be the global multiscale basis functions obtained in (\ref{msbasis2}). Then we have
\[
\norm{\psi^i_j-\psi^i_{j,l}}_a^2+\norm{\pi(\psi^i_j-\psi^i_{j,l})}_s^2\leq \frac{15}{4}(1+\frac{1}{\Lambda})\theta^l\norm{\phi^i_j}_s^2,
\]
where $\theta=\tilde{c}/(1+\tilde{c})$ and $\tilde{c}=2(1+1/\sqrt{\Lambda})$.
\end{theorem}
\begin{proof}
In terms of (\ref{msbasis1}) and (\ref{msbasis2}), we have (let the oversampling domain considered in the variational problem (\ref{msbasis1}) be $K_j^{l+1}$, i.e., find $\psi^i_{j,l}\in V(K_j^{l+1})$ such that $a(\psi_{j,l}^i,v)+s(\pi\psi_{j,l}^i,\pi v)=s(\phi_j^i,\pi v),$ for all $v\in V(K_j^{l+1})$)
\[
a(\psi^i_j-\psi_{j,l}^i,v)+s(\pi(\psi^i_j-\psi_{j,l}^i),\pi v)=0,\quad\forall v\in V(K_j^{l+1}).
\]
Taking $v=w-\psi^i_{j,l}$ with $w\in V(K_j^{l+1})$ in the above relation, we have 
\begin{equation*}
	\begin{aligned}
		&\quad\norm{\psi_j^i-\psi_{j,l}^i}_a^2+\norm{\pi(\psi_j^i-\psi_{j,l}^i)}_s^2\\
  &=a(\psi_j^i-\psi_{j,l}^i,\psi_j^i-w+w-\psi_{j,l}^i)+s(\pi(\psi_j^i-\psi_{j,l}^i),\pi(\psi_j^i-w+w-\psi_{j,l}^i))\\
  &=a(\psi_j^i-\psi_{j,l}^i,\psi_j^i-w)+s(\pi(\psi_j^i-\psi_{j,l}^i),\pi(\psi_j^i-w))\\
  &\leq (\norm{\psi_j^i-\psi_{j,l}^i}_a^2+\norm{\psi_j^i-w}_a^2)/2+(\norm{\pi(\psi_j^i-\psi_{j,l}^i)}_s^2+\norm{\pi(\psi_j^i-w)}_s^2)/2.
	\end{aligned}
\end{equation*}
Thus, $\quad\norm{\psi_j^i-\psi_{j,l}^i}_a^2+\norm{\pi(\psi_j^i-\psi_{j,l}^i)}_s^2\leq \norm{\psi_j^i-w}_a^2+\norm{\pi(\psi_j^i-w)}_s^2$
for all $w\in V(K_j^{l+1})$. Letting $w=\eta_j^{l+1,l}\psi^i_j$ in the above relation, we have
\begin{equation}
\label{estimate_1}
\quad\norm{\psi_j^i-\psi_{j,l}^i}_a^2+\norm{\pi(\psi_j^i-\psi_{j,l}^i)}_s^2\leq \norm{\psi_j^i-\eta_j^{l+1,l}\psi^i_j}_a^2+\norm{\pi(\psi_j^i-\eta_j^{l+1,l}\psi^i_j)}_s^2.
\end{equation}
Next we will estimate the two terms on the right hand side of (\ref{estimate_1}). We divide the proof into five steps.

\textbf{Step 1:} We will estimate the first term on the right hand side of (\ref{estimate_1}). By the definition of the norm $\norm{\cdot}_a$ and the fact that $\text{supp}(1-\eta_j^{l+1,l})\subset(\Omega\backslash K_j^l)$, we have
\begin{equation}
\label{1-eta times psi}
\begin{aligned}
\norm{(1-\eta_j^{l+1,l})\psi^i_j}_a^2&\leq\int_{\Omega\backslash K_j^l}\varepsilon^2\abs{1-\eta_j^{l+1,l}}^2\cdot\abs{\nabla\psi^i_j}^2dx+\int_{\Omega\backslash K_j^l}\varepsilon^2\abs{\nabla\eta_j^{l+1,l}}^2\cdot\abs{\psi^i_j}^2dx\\
&\quad+\int_{\Omega\backslash K_j^l}\varepsilon^2(1-\eta_j^{l+1,l})^2\cdot\abs{\psi^i_j}^2dx\\
&\leq 2(\norm{\psi}_{a(\Omega\backslash K_j^l)}^2+\norm{\psi}_{s(\Omega\backslash K_j^l)}^2)
\end{aligned}
\end{equation}
since $\frac{1}{2}\varepsilon^2\abs{\nabla\eta_j^{l+1,l}}^2\leq \tilde{\mu}(x)$ by the construction of $\eta_j^{l+1,l}$. Note that for each $K_j$ ($1\leq j\leq N$), by Lemma \ref{inter}, we have
\begin{equation}
\begin{aligned}
\label{psi_s_k}
\norm{\psi_j^i}_{s(K_j)}^2&=\norm{(I-\pi)\psi^i_j+\pi(\psi^i_j)}_{s(K_j)}^2\\
&=\norm{(I-\pi)\psi^i_j}_{s(K_j)}^2+\norm{\pi(\psi^i_j)}_{s(K_j)}^2\\
&\leq \frac{1}{\Lambda}\norm{\psi^i_j}_{a(K_j)}^2+\norm{\pi(\psi^i_j)}_{s(K_j)}^2.
\end{aligned}
\end{equation}
Summing (\ref{psi_s_k}) over all $K_j\subset \Omega\backslash K_j^l$ and combining (\ref{1-eta times psi}), we obtain that
\begin{align}
\label{estimat for the first term before Step 1}
\norm{(1-\eta_j^{l+1,l})\psi^i_j}_a^2\leq 2(1+\frac{1}{\Lambda})(\norm{\psi^i_j}_{a(\Omega\backslash K_j^l)}^2+\norm{\pi(\psi^i_j)}_{s(\Omega\backslash K_j^l)}^2).
\end{align}

\textbf{Step 2:} We will estimate the second term on the right hand side of (\ref{estimate_1}).  By Lemma \ref{inter}, we know 
\[
\norm{\pi(\psi_j^i-\eta_j^{l+1,l}\psi^i_j)}_s^2\leq \norm{(1-\eta_j^{l+1,l})\psi^i_j}_s^2\leq \norm{\psi^i_j}_{s(\Omega\backslash K_j^l)}^2.
\]
By using (\ref{psi_s_k}), we have
\begin{equation}
\label{estimat for the first term before Step 2}
\norm{\pi(\psi_j^i-\eta_j^{l+1,l}\psi^i_j)}_s^2\leq\frac{1}{\Lambda}\norm{\psi^i_j}_{a(\Omega\backslash K_j^l)}^2+\norm{\pi(\psi^i_j)}_{s(\Omega\backslash K_j^l)}^2.
\end{equation}

In terms of Steps 1 and 2 (i.e. by (\ref{estimat for the first term before Step 1}) and (\ref{estimat for the first term before Step 2})), we see that (\ref{estimate_1}) can be estimated as
\[
\quad\norm{\psi_j^i-\psi_{j,l}^i}_a^2+\norm{\pi(\psi_j^i-\psi_{j,l}^i)}_s^2\leq 3(1+\frac{1}{\Lambda})(\norm{\psi^i_j}_{a(\Omega\backslash K_j^l)}^2+\norm{\pi(\psi^i_j)}_{s(\Omega\backslash K_j^l)}^2).
\]

\textbf{Step 3:} In this step, we estimate $\norm{\psi^i_j}_{a(\Omega\backslash K_j^l)}^2+\norm{\pi(\psi^i_j)}_{s(\Omega\backslash K_j^l)}^2$. By utilizing (\ref{msbasis2}) and the test function $(1-\eta_j^{l,l-1})\psi^i_j$, we have
\begin{equation}
\label{use the test function to global basis}
\begin{aligned}
&\quad a(\psi_{j}^i,(1-\eta_j^{l,l-1})\psi^i_j)+s(\pi\psi_{j}^i,\pi ((1-\eta_j^{l,l-1})\psi^i_j))\\
&=s(\phi_{j}^i,\pi ((1-\eta_j^{l,l-1})\psi^i_j))=0
\end{aligned}
\end{equation}

since $\text{supp}(1-\eta_j^{l,l-1})\subset(\Omega\backslash K_j^{l-1})$ and $\text{supp}(\phi^i_j)\subset K_j$. Notice that
\[
a(\psi_{j}^i,(1-\eta_j^{l,l-1})\psi^i_j)=\int_{\Omega\backslash K_j^{l-1}}\frac{1}{2}\varepsilon^2\nabla\psi^i_j\cdot\nabla(\overline{(1-\eta_j^{l,l-1})\psi^i_j})+V^\delta(x)\psi^i_j\cdot(1-\eta_j^{l,l-1})\overline{\psi^i_j},
\]
so we have
\begin{align*}
a(\psi_{j}^i,(1-\eta_j^{l,l-1})\psi^i_j)&=\int_{\Omega\backslash K_j^{l-1}}\frac{1}{2}\varepsilon^2\abs{\nabla\psi^i_j}^2\cdot(1-\eta_j^{l,l-1})-\int_{\Omega\backslash K_j^{l-1}}\frac{1}{2}\varepsilon^2\nabla\psi^i_j\cdot\nabla\eta_j^{l,l-1}\cdot\overline{\psi^i_j}\\
&\quad+\int_{\Omega\backslash K_j^{l-1}}V^\delta(x)\abs{\psi^i_j}^2\cdot(1-\eta_j^{l,l-1}).
\end{align*}
Then we obtain that 
\begin{align*}
\norm{\psi^i_j}_{a(\Omega\backslash K_j^l)}^2&\leq\int_{\Omega\backslash K_j^{l-1}}\frac{1}{2}\varepsilon^2(1-\eta_j^{l,l-1})\abs{\nabla\psi^i_j}^2+\int_{\Omega\backslash K_j^{l-1}}V^\delta(x)(1-\eta_j^{l,l-1})\abs{\psi^i_j}^2\\
&=a(\psi_{j}^i,(1-\eta_j^{l,l-1})\psi^i_j)+\int_{\Omega\backslash K_j^{l-1}}\frac{1}{2}\varepsilon^2\nabla\psi^i_j\cdot\nabla\eta_j^{l,l-1}\cdot\overline{\psi^i_j}\\
&\leq a(\psi_{j}^i,(1-\eta_j^{l,l-1})\psi^i_j)+\norm{\psi^i_j}_{a(K_j^l\backslash K_j^{l-1})}\cdot\norm{\psi^i_j}_{s(K_j^l\backslash K_j^{l-1})}.
\end{align*}
Since $\eta_j^{l,l-1}\equiv0$ in $\Omega\backslash K_j^l$, we have
\[
s(\pi(\psi^i_j),\pi((1-\eta_j^{l,l-1})\psi^i_j))=\norm{\pi(\psi^i_j)}_{s(\Omega\backslash K_j^l)}^2+\int_{K_j^l\backslash K_j^{l-1}}\tilde{\mu}(x)\pi(\psi^i_j)\cdot\overline{\pi((1-\eta_j^{l,l-1})\psi^i_j)}.
\]
Thus, combining Lemma \ref{inter}, we obtain
\begin{align*}
\norm{\pi(\psi^i_j)}_{s(\Omega\backslash K_j^l)}^2&=s(\pi(\psi^i_j),\pi((1-\eta_j^{l,l-1})\psi^i_j))-\int_{K_j^l\backslash K_j^{l-1}}\tilde{\mu}(x)\pi(\psi^i_j)\cdot\overline{\pi((1-\eta_j^{l,l-1})\psi^i_j)}\\
&\leq s(\pi(\psi^i_j),\pi((1-\eta_j^{l,l-1})\psi^i_j))+\norm{\pi(\psi^i_j)}_{s(K_j^l\backslash K_j^{l-1})}\cdot\norm{\psi^i_j}_{s(K_j^l\backslash K_j^{l-1})}.
\end{align*}
{Summing the above} two inequalities and in terms of (\ref{psi_s_k}) , we have
\begin{align*}
&\quad\norm{\psi^i_j}_{a(\Omega\backslash K_j^l)}^2+\norm{\pi(\psi^i_j)}_{s(\Omega\backslash K_j^l)}^2\\
&\leq \norm{\psi^i_j}_{s(K_j^l\backslash K_j^{l-1})}(\norm{\psi^i_j}_{a(K_j^l\backslash K_j^{l-1})}+\norm{\pi(\psi^i_j)}_{s(K_j^l\backslash K_j^{l-1})})\\
&\leq(\frac{1}{\sqrt{\Lambda}}\norm{\psi^i_j}_{a(K_j^l\backslash K_j^{l-1})}+\norm{\pi(\psi^i_j)}_{s(K_j^l\backslash K_j^{l-1})})\cdot(\norm{\psi^i_j}_{a(K_j^l\backslash K_j^{l-1})}+\norm{\pi(\psi^i_j)}_{s(K_j^l\backslash K_j^{l-1})})\\
&\leq 2(\frac{1}{\sqrt{\Lambda}}+1)(\norm{\psi^i_j}_{a(K_j^l\backslash K_j^{l-1})}^2+\norm{\pi(\psi^i_j)}_{s(K_j^l\backslash K_j^{l-1})}^2).
\end{align*}

\textbf{Step 4:} In this step, we will show that $\norm{\psi^i_j}_{a(\Omega\backslash K_j^l)}^2+\norm{\pi(\psi^i_j)}_{s(\Omega\backslash K_j^l)}^2$ can be estimated by $\norm{\psi^i_j}_{a(\Omega\backslash K_j^{l-1})}^2+\norm{\pi(\psi^i_j)}_{s(\Omega\backslash K_j^{l-1})}^2$. Based on Step 3, we have
\begin{align*}
&\quad \norm{\psi^i_j}_{a(\Omega\backslash K_j^{l-1})}^2+\norm{\pi(\psi^i_j)}_{s(\Omega\backslash K_j^{l-1})}^2\\
&=\norm{\psi^i_j}_{a(\Omega\backslash K_j^l)}^2+\norm{\pi(\psi^i_j)}_{s(\Omega\backslash K_j^l)}^2+\norm{\psi^i_j}_{a(K_j^l\backslash K_j^{l-1})}^2+\norm{\pi(\psi^i_j)}_{s(K_j^l\backslash K_j^{l-1})}^2\\
&\geq (1+(2(\frac{1}{\sqrt{\Lambda}}+1))^{-1})(\norm{\psi^i_j}_{a(\Omega\backslash K_j^l)}^2+\norm{\pi(\psi^i_j)}_{s(\Omega\backslash K_j^l)}^2).
\end{align*}
Utilizing the above inequality recursively, we have
\[
\norm{\psi^i_j}_{a(\Omega\backslash K_j^l)}^2+\norm{\pi(\psi^i_j)}_{s(\Omega\backslash K_j^l)}^2\leq (1+(2(\frac{1}{\sqrt{\Lambda}}+1))^{-1})^{-l}(\norm{\psi^i_j}_a^2+\norm{\pi(\psi^i_j)}_s^2).
\]
Combining the above steps, we get
\[
\norm{\psi_j^i-\psi_{j,l}^i}_a^2+\norm{\pi(\psi_j^i-\psi_{j,l}^i)}_s^2\leq 3(1+\frac{1}{\Lambda})\theta^l(\norm{\psi^i_j}_a^2+\norm{\pi(\psi^i_j)}_s^2),
\]
where $\theta=\tilde{c}/(1+\tilde{c})$, and $\tilde{c}=2(\frac{1}{\sqrt{\Lambda}}+1)$. Note that $0<\theta<1$.

\textbf{Step 5:} In this step, we prove the final result. By (\ref{msbasis2}), we have 
\begin{equation*}
a(\psi_{j}^{i},\psi_{j}^{i})+s\big(\pi(\psi_{j}^{i}),\pi(\psi_{j}^{i})\big)=s\big(\phi^i_j,\pi(\psi_{j}^{i})\big)
\end{equation*}
By the above equality, on one hand, we have
$\norm{\pi(\psi_{j}^i)}_s^2\leq\norm{\phi^i_j}_s\cdot\norm{\pi(\psi_{j}^i)}_s$. That is,
\begin{equation}
\label{con_11}
\norm{\pi(\psi_{j}^i)}_s^2\leq\norm{\phi^i_j}_s^2.
\end{equation}
On the other hand, 
\begin{align*}
s\big(\phi^i_j,\pi(\psi_{j}^{i})\big)&\leq\frac{1}{\sqrt{2}}\norm{\phi^i_j}_s\cdot\sqrt{2}\norm{\pi(\psi_{j}^i)}_s\leq(\frac{1}{2}\norm{\phi^i_j}_s^2+2\norm{\pi(\psi_{j}^i)}_s^2)/2\\
&=\frac{1}{4}\norm{\phi^i_j}_s^2+s\big(\pi(\psi_{j}^{i}),\pi(\psi_{j}^{i})\big)
\end{align*}
Then, we have
\begin{equation}
\label{con_22}
\norm{\psi_{j}^i}_a^2\leq\frac{1}{4}\norm{\phi^i_j}_s^2.
\end{equation}
By combining (\ref{con_11}), (\ref{con_22}) and Step 4, we obtain
\begin{equation*}
\norm{\psi_j^i-\psi_{j,l}^i}_a^2+\norm{\pi(\psi_j^i-\psi_{j,l}^i)}_s^2\leq \frac{15}{4}(1+\frac{1}{\Lambda})\theta^l\norm{\phi_j^i}_s^2,
\end{equation*}
where $\theta=\tilde{c}/(1+\tilde{c})$, and $\tilde{c}=2(\frac{1}{\sqrt{\Lambda}}+1)$.
\end{proof}

Next, we prove the convergence of the multiscale solution $u_\text{ms}$ to $u$. Before this, we give an assumption.

\begin{assumption}
\label{oversampling layers}
The oversampling size $l$ satisfies
\[
(l+1)^{\frac{d}{2}}\cdot \theta^{\frac{l}{2}}\cdot\frac{1}{V_\text{min}}=O(\frac{H^2}{\varepsilon^2}),
\]
where $\theta=\tilde{c}/(1+\tilde{c})$, and $\tilde{c}=2(\frac{1}{\sqrt{\Lambda}}+1)$.
\end{assumption}

\begin{theorem}
\label{convergence of the multiscale solution}
Let $u$ be the solution of the stationary problem (\ref{stationary problem}) and $u_\text{ms}$ be the multiscale solution of (\ref{multiscale problem}). Under Assumptions \ref{regularity for mesh} and \ref{oversampling layers}, we have
\begin{equation}
\label{energy norm for u-u_ms}
\norm{u-u_\text{ms}}_a\leq C/\sqrt{\Lambda}\norm{f}_{s^{-1}},
\end{equation}
and
\begin{equation}
\label{s norm for u-u_ms}
\norm{u-u_\text{ms}}_s\leq C/\Lambda\norm{f}_{s^{-1}},
\end{equation}
where $C$ is independent of  $\Lambda,\varepsilon,\delta,H$, and the overlampling size $l$.
\end{theorem}
\begin{proof}
In terms of Lemma \ref{decay property} and the fact that $s(\phi^i_j,\phi^i_j)=1$, we have
\[
\norm{\psi^i_j-\psi^i_{j,l}}_a^2+\norm{\pi(\psi^i_j-\psi^i_{j,l})}_s^2\leq \frac{15}{4}(1+\frac{1}{\Lambda})\theta^l,
\]
where $\theta=\tilde{c}/(1+\tilde{c})$ and $\tilde{c}=2(1+1/\sqrt{\Lambda})$. We write $u_\text{glo}=\sum_{j=1}^N\sum_{i=1}^{l_j}c^i_j\psi_j^i$ and define $w=\sum_{j=1}^N\sum_{i=1}^{l_j}c^i_j\psi_{j,l}^i\in V_\text{ms}$. By the orthogonal property, we have (Note that $w-u_\text{ms}\in V_\text{ms}$)
\begin{align*}
\norm{u-u_\text{ms}}_a^2&=a(u-u_\text{ms},u-u_\text{ms})=a(u-u_\text{ms},u-w+w-u_\text{ms})\\
&=a(u-u_\text{ms},u-w)\leq \norm{u-u_\text{ms}}_a\cdot\norm{u-w}_a.
\end{align*}
Then we have
$$\norm{u-u_\text{ms}}_a\leq\norm{u-w}_a\leq\norm{u-u_\text{glo}}_a+\norm{\sum_{i=1}^N\sum_{j=1}^{l_i}c^i_j(\psi_j^i-\psi_{j,l}^i)}_a.$$
By utilizing Assumption \ref{regularity for mesh}, Lemma \ref{decay property}, the fact that $s(\phi^i_j,\phi^l_j)=\delta_{il}$,  $\big($applying them to the function $\sum_{j=1}^{l_i}c^i_j(\psi_j^i-\psi_{j,l}^i)$$\big)$ and denoting $\phi=\sum_{i=1}^N\sum_{j=1}^{l_i}c^i_j\phi^i_j\in V_\text{aux}$, we obtain
\begin{align*}
\norm{w-u_\text{glo}}_a^2&\leq C_\text{ol}(l+1)^d\sum_{i=1}^N\norm{\sum_{j=1}^{l_i}c^i_j(\psi_j^i-\psi_{j,l}^i)}_a^2\\
&\leq C_\text{ol}(l+1)^d\cdot\frac{15}{4}(1+\frac{1}{\Lambda})\theta^l\sum_{i=1}^N\sum_{j=1}^{l_i}\norm{c^i_j\phi_j^{i}}_s^2\\
&=C_\text{ol}(l+1)^d\cdot\frac{15}{4}(1+\frac{1}{\Lambda})\theta^l\sum_{i=1}^N\sum_{j=1}^{l_i}(c^i_j)^2\\
&=C_\text{ol}(l+1)^d\cdot\frac{15}{4}(1+\frac{1}{\Lambda})\theta^l\cdot s(\phi,\phi).
\end{align*}
By the definitions of $u_\text{glo},\phi$ and the variational form (\ref{msbasis2}), we know 
\begin{equation}
\label{u glo variational}
a(u_\text{glo},v)+s(\pi(u_\text{glo}),\pi(v))=s(\phi,\pi(v)),\quad \forall v\in V.
\end{equation}
For $\phi\in V_\text{aux}$, by Lemma \ref{C inverse}, there is $\alpha\in V$ such that
$\pi(\alpha)=\phi$, $\norm{\alpha}_a^2\leq C_0\norm{\phi}_s^2$.
Letting $v=\alpha$ in (\ref{u glo variational}), we have
$
a(u_\text{glo},\alpha)+s(\pi(u_\text{glo}),\pi(\alpha))=s(\phi,\pi(\alpha))=s(\phi,\phi)
$.
Then we obtain that
\begin{align*}
s(\phi,\phi)&=a(u_\text{glo},\alpha)+s(\pi(u_\text{glo}),\pi(\alpha))\leq\norm{u_\text{glo}}_a\norm{\alpha}_a+\norm{\pi(u_\text{glo})}_s\norm{\phi}_s\\
&\leq(\sqrt{C_0}+1)\norm{\phi}_s(\norm{u_\text{glo}}_a+\norm{\pi(u_\text{glo})}_s).
\end{align*}
Therefore, we have
$$\norm{w-u_\text{glo}}_a\leq \frac{\sqrt{15}}{2}C_\text{ol}^{\frac{1}{2}}(l+1)^{d/2}(\sqrt{C_0}+1)(1+\frac{1}{\Lambda})^{\frac{1}{2}}\theta^{l/2}(\norm{u_\text{glo}}_a+\norm{\pi(u_\text{glo})}_s).$$ Next we estimate $\norm{u_\text{glo}}_a+\norm{\pi(u_\text{glo})}_s$. In terms of Lemma \ref{inter} and (\ref{global multiscale problem}), we have
\begin{align*}
\norm{\pi(u_\text{glo})}_s^2&\leq\norm{u_\text{glo}}_s^2=s(u_\text{glo},u_\text{glo})\leq \max\{\tilde{\mu}(x)\}\norm{u_\text{glo}}_{L^2}^2\\
&\leq \max\{\tilde{\mu}(x)\}\cdot \frac{1}{V_\text{min}}\norm{u_\text{glo}}_a^2=\max\{\tilde{\mu}(x)\}\cdot \frac{1}{V_\text{min}}(f,u_\text{glo})\\
&\leq \max\{\tilde{\mu}(x)\}\cdot \frac{1}{V_\text{min}}\norm{u_\text{glo}}_s\cdot\norm{f}_{s^{-1}}.
\end{align*}
So we obtain that
\[
\norm{\pi(u_\text{glo})}_s\leq\norm{u_\text{glo}}_s\leq \max\{\tilde{\mu}(x)\}\cdot \frac{1}{V_\text{min}}\norm{f}_{s^{-1}}.
\]
For $\norm{u_\text{glo}}_a$, we have
\[
\norm{u_\text{glo}}_a^2=a(u_\text{glo},u_\text{glo})=(f,u_\text{glo})\leq\norm{f}_{s^{-1}}\cdot\norm{u_\text{glo}}_s\leq \max\{\tilde{\mu}(x)\}\cdot \frac{1}{V_\text{min}}\norm{f}_{s^{-1}}^2.
\]
Then we obtain that (Note that we assume $\varepsilon/H$ can be controlled by a constant)
\[
\norm{u_\text{glo}}_a+\norm{\pi(u_\text{glo})}_s\leq C^*\max\{\tilde{\mu}(x)\}\cdot \frac{1}{V_\text{min}}\norm{f}_{s^{-1}},
\]
where $C^*$ in independent of $\varepsilon$, $H$. Therefore, combining Theorem \ref{glothm}, we get
\[
\norm{u-u_\text{ms}}_a\leq (\frac{1}{\sqrt{\Lambda}}+\frac{\sqrt{15}}{2}C_\text{ol}^{\frac{1}{2}}(l+1)^{d/2}(\sqrt{C_0}+1)(1+\frac{1}{\Lambda})^{\frac{1}{2}}\theta^{l/2}C^*\max\{\tilde{\mu}(x)\}\cdot \frac{1}{V_\text{min}})\norm{f}_{s^{-1}}.
\]
Recalling Assumption \ref{oversampling layers}, we have $$C_\text{ol}^{\frac{1}{2}}(l+1)^{d/2}(\sqrt{C_0}+1)(1+\frac{1}{\Lambda})^{\frac{1}{2}}\theta^{l/2}C^*\max\{\tilde{\mu}(x)\}\cdot \frac{1}{V_\text{min}}\sim O(1)$$
since $\max\{\tilde{\mu}(x)\}=O(\varepsilon^2/H^2)$ and $C_\text{ol},C_0,C^*$ are $O(1)$. Thus we finally obtain that $\norm{u-u_\text{ms}}_a\leq C\Lambda^{-\frac{1}{2}}\norm{f}_{s^{-1}},$
where $C$ does not depend on $\varepsilon,\delta,H,\Lambda$. 

For the estimate $\norm{u-u_\text{ms}}_s\leq C/\Lambda\norm{f}_{s^{-1}}$, a similar technique to Theorem \ref{glothm} can be applied.
This completes the proof.
\end{proof}
\section{{Convergence of Crank-Nicolson CEM scheme}}\label{sec5}   
In this section, we analyze the convergence of the constraint energy minimizing generalized multiscale finite element method for (\ref{model problem}). We assume Assumptions \ref{assumption for u0 and potential}-\ref{oversampling layers} are satisfied. 
\subsection{Regularity of the solution} %to the Schr\"{o}dinger equation with multiscale potential}
Similar to \cite{wu2022}, we study the temporal regularity of the solution $u$ of the Schr\"{o}dinger equation (\ref{model problem}). The spatial regularity of the solution $u$ is analogous so we omit details for the temporal regularity.

\begin{lemma}
\label{temporal regularity}
Let $u$ be the solution of (\ref{model problem}). If $\partial^k_tu(t)\in L^2(\Omega)$ {($L^2(\Omega)$ is the standard Lebesgue space)} for any $t\in[0,T]$, where $k=1,2,3$, then it holds true that for any $0\leq t\leq T$
\[
\norm{\partial^k_tu(t)}_{s^{-1}}\leq \frac{CH\varepsilon^{k-3}}{\min\{\varepsilon^{2k-2},\delta^{2k-2}\}},
\]
{where $\norm{v}_{s^{-1}}=O(H/\varepsilon)$ for all $v\in L^2(\Omega)$ by the definition of $s^{-1}$-norm (i.e. the scaling term $H/\varepsilon$ is hidden in the norm $\norm{\cdot}_{s^{-1}}$ for all functions in the standard Lebesgue space $L^2(\Omega)$).}
\end{lemma}

\begin{proof}
The proof is analogous to \cite[Lemma 4.1]{wu2022}.  In terms of (\ref{model problem}), we have
\[
i\varepsilon\int_{\Omega}\partial^2_tu\cdot\overline{\partial_tu}dx=\frac{1}{2}\varepsilon^2\int_{\Omega}\nabla (\partial_tu)\cdot\nabla(\overline{\partial_tu})dx+\int_{\Omega}V^\delta(x)\partial_tu\cdot\overline{\partial_tu}dx.
\]
Taking the imaginary part, we have $\frac{d}{dt}\norm{\partial_tu}^2=0$, which implies $\norm{\partial_tu(t)}=\norm{\partial_tu(0)}$, for all $t\in[0,T]$. Since $i\varepsilon\partial_tu(0)=-\frac{\varepsilon^2}{2}\Delta u_0+Vu_0$ and $\norm{\cdot}_{s^{-1}}= O(H/\varepsilon)$, then we have
\[
\norm{\partial_tu(0)}_{s^{-1}}\leq \frac{\varepsilon}{2}\norm{\Delta u_0}_{s^{-1}}+\frac{1}{\varepsilon}\norm{Vu_0}_{s^{-1}}\leq \frac{CH}{\varepsilon^2}.
\]
Doing similar process to $\partial^k_tu$ for $k=2,3,4$ and using Assumption \ref{assumption for u0 and potential}, we can obtain the desired results.
\end{proof}
\subsection{Projection error}
Let $u(t)$ be the solution of Schr\"{o}dinger equation (\ref{model problem}) and $\sigma(u(t))$ be the projection of $u(t)$ in $V_\text{ms}$ such that for all $0\leq t\leq T$,
\begin{equation}
\label{projection to multiscale space}
a(u(t)-\sigma(u(t)),w)=0,\quad \forall w\in V_\text{ms}.
\end{equation}
Next we give the following lemma on the projection error estimates.
\begin{lemma}
\label{regularity for projection operator}
If $\partial^{k+1}_tu(t)\in L^2(\Omega)$ ($k=0,1,2$,  {$L^2(\Omega)$ is the standard Lebesgue space}) for any $t\in[0,T]$, then it holds true that for any $0\leq t\leq T$
\[
\norm{\partial^k_tu(t)-\partial^k_t\sigma(u(t))}_a\leq \frac{CH}{\sqrt{\Lambda}\varepsilon^{1-k}\min\{\varepsilon^{2k},\delta^{2k}\}},
\]
and
\begin{equation*}
\label{s norm for u-u_pro}
\norm{\partial^k_tu(t)-\partial^k_t\sigma(u(t))}_s\leq \frac{CH}{\Lambda \varepsilon^{1-k}\min\{\varepsilon^{2k},\delta^{2k}\}},
\end{equation*}
where $C$ is independent of  $\Lambda,\varepsilon,\delta,H$, and the {overlampling size $l$.} {$\norm{v}_s=O(\varepsilon/H)$ for all $v\in L^2(\Omega)$ by the definition of $s$-norm (i.e. the scaling term $\varepsilon/H$ is hidden in the norm $\norm{\cdot}_s$ for all functions in the standard Lebesgue space $L^2(\Omega)$).}
\end{lemma}
\begin{proof}
For simplicity, we just consider the case where $k=0$ (i.e. $\partial_tu(t)\in L^2(\Omega)$, $\forall t\in[0,T]$). The proof of the case where $k=1,2$ is simlar to that of $k=0$. In terms of Theorem \ref{convergence of the multiscale solution}, Lemma \ref{temporal regularity} and (\ref{projection to multiscale space}), we have that for any $0\leq t\leq T$,
\begin{align*}
&\quad\norm{u(t)-\sigma(u(t))}_a^2=a(u(t)-\sigma(u(t)),u(t)-\sigma(u(t)))\\
&=a(u(t)-u_\text{ms}(t),u(t)-\sigma(u(t)))+a(u_\text{ms}(t)-\sigma(u(t)),u(t)-\sigma(u(t)))\\
&=a(u(t)-u_\text{ms}(t),u(t)-\sigma(u(t)))\leq \norm{u(t)-u_\text{ms}(t)}_a\cdot\norm{u(t)-\sigma(u(t))}_a.
\end{align*}
Then we obtain that
\begin{align*}
\norm{u(t)-\sigma(u(t))}_a&\leq\norm{u(t)-u_\text{ms}(t)}_a\leq C\frac{1}{\sqrt{\Lambda}}\norm{\mathcal{H}(u(t))}_{s^{-1}}\\
&\leq C\frac{\varepsilon}{\sqrt{\Lambda}}\norm{\partial_tu(t)}_{s^{-1}}\leq C\frac{\varepsilon}{\sqrt{\Lambda}}\cdot\frac{H}{\varepsilon^2}=C\Lambda^{-\frac{1}{2}}\frac{H}{\varepsilon}.
\end{align*}
For the estimate $\norm{\partial^k_tu(t)-\partial^k_t\sigma(u(t))}_s$, a similar technique to Theorem \ref{glothm} can be applied. This completes the proof.
\end{proof}
\subsection{Crank-Nicolson CEM scheme}
We analyze the error estimate of the Crank-Nicolson CEM scheme for Schr\"{o}dinger equation (\ref{model problem}). More precisely, we utilize the Crank-Nicolson scheme in temporal discretization and the constraint energy minimizing generalized multiscale finite element method in spatial discretization. For some $N\in\mathbb{N}$ and $N>0$, let $\Delta t=T/N$ and $t_n=n\Delta t$, $n=0,1,\dots, N$. For the multiscale space $V_\text{ms}$, we approximate $u^n\in V_{\text{ms}}$ such that 
\begin{equation}
\label{C-N1}
\begin{cases}
\mathrm{i}\varepsilon(\frac{u^n-u^{n-1}}{\Delta t},\overline{v})=\frac{1}{2}\varepsilon^2(\frac{\nabla u^n+\nabla u^{n-1}}{2},\overline{\nabla v})+(\frac{V^\delta u^n+V^\delta u^{n-1}}{2},\overline{v}),\\
u^0=\sigma(u(t_0))
\end{cases}
\end{equation}
for  all $v\in V_{\text{ms}}$ and $n=1,2,...,N$. {For the Crank-Nicolson CEM scheme (\ref{C-N1}), we emphasize that the ``CEM" designation refers to the constraint energy minimizing generalized multiscale finite element method used for constructing the spatial approximation space $V_{\text{ms}}$, as detailed in Section \ref{sec3}.}

\begin{theorem}
\label{final convergence}
Let $u^n$ be the solution of  (\ref{C-N1}) and $u$ the solution of (\ref{model problem}). If $\partial_t^ku(t)\in L^2(\Omega)$ ($k=1,2,3$) for any $t\in[0,T]$, then we have that
\begin{equation}
\label{L2 error}
\norm{u^N-u(T)}\leq C(\frac{H^2}{\Lambda\sqrt{\varepsilon}\min\{\varepsilon^2,\delta^2\}}+\frac{\varepsilon^{3/2}\Delta t^2}{\min\{\varepsilon^4,\delta^4\}}).
\end{equation}
In addition, we have $\int_0^T\norm{e(t)}_a\leq C(\frac{TH^2}{\Lambda\sqrt{\varepsilon}\min\{\varepsilon^2,\delta^2\}}+\frac{TH}{\sqrt{\Lambda\varepsilon}}+\frac{T\varepsilon^{3/2} \Delta t^2}{\min\{\varepsilon^4,\delta^4\}})$, where $e(t)=u^n-u(t_n)$ for $t\in(t_{n-1},t_n]$.
Here $C$ is independent of  $\Lambda,\varepsilon,\delta,H$, and the {oversampling size $l$.}
\end{theorem}
\begin{remark}
\label{relations between parameters}
We set $\beta=\frac{H}{\sqrt{\Lambda}}$, allowing us to express the related error terms as $\beta^2\varepsilon^{-\frac{5}{2}}+\beta\varepsilon^{-1/2}$ if $\varepsilon\leq\delta$ (or $\beta^2\varepsilon^{-\frac{1}{2}}\delta^{-2}+\beta\varepsilon^{-1/2}$ if $\delta<\varepsilon$). By choosing $\beta$ sufficiently small, we ensure that $\beta^2\varepsilon^{-\frac{5}{2}}+\beta\varepsilon^{-1/2}$ (or $\beta^2\varepsilon^{-\frac{1}{2}}\delta^{-2}+\beta\varepsilon^{-1/2}$) is sufficiently small. On the other hand, one can select $\Delta t$ small enough to guarantee $\frac{\varepsilon^{3/2}\Delta t^2}{\min\{\varepsilon^4,\delta^4\}}$ remaining sufficiently small. {More precisely, the convergence requires $H/\sqrt{\Lambda}=O(\varepsilon^{\frac{5}{4}})$, $\Delta t=O(\varepsilon^{\frac{5}{4}})$ if $\varepsilon\leq \delta$; while if $\delta<\varepsilon$, the convergence requires $H/\sqrt{\Lambda}=O(\varepsilon^{\frac{1}{4}}\delta)$, $\Delta t=O(\frac{\delta^2}{\varepsilon^{3/4}})$.} % Thus we conclude that the error for $u^n-u(T)$ is sufficiently small. We remark that in our numerical experiments, we mainly consider moderately small Planck constant $\varepsilon$, but not too small. 
\end{remark}
\begin{proof}[proof of Theorem \ref{final convergence}]
In terms of (\ref{model problem}) and (\ref{projection to multiscale space}), we have
\begin{equation}
\label{projection equality and model problem}
a(\frac{\sigma(u(t_n))+\sigma(u(t_{n-1}))}{2},\overline{v})=a(\frac{u(t_n)+u(t_{n-1})}{2},\overline{v})=i\varepsilon(\frac{\partial_tu(t_n)+\partial_tu(t_{n-1})}{2},v).
\end{equation}
Letting $e^n=u^n-\sigma(u(t_n))$ and subtracting (\ref{projection equality and model problem}) from (\ref{C-N1}), we have
\begin{equation}
\label{en form}
a(\frac{e^n+e^{n-1}}{2},\overline{v})=i\varepsilon(\frac{e^n-e^{n-1}}{\Delta t},\overline{v})+i\varepsilon(y_1^n+y_2^n,\overline{v}),
\end{equation}
where $y^n_1=\frac{\sigma(u(t_n))-\sigma(u(t_{n-1}))}{\Delta t}-\frac{u(t_n)-u(t_{n-1})}{\Delta t}$, $y^n_2=\frac{u(t_n)-u(t_{n-1})}{\Delta t}-\frac{\partial_tu(t_n)+\partial_tu(t_{n-1})}{2}$. Let $v=e^n+e^{n-1}$, then we have
\begin{equation}
\label{enplus test}
\frac{1}{2}\norm{e^n+e^{n-1}}_a^2=i\varepsilon\cdot\frac{1}{\Delta t}(\norm{e^n}^2-\norm{e^{n-1}}^2)+i\varepsilon(y_1^n+y_2^n,\overline{e^n+e^{n-1}}).
\end{equation}
Taking imaginary part of (\ref{enplus test}), we get
$\frac{1}{\Delta t}(\norm{e^n}^2-\norm{e^{n-1}}^2)\leq \norm{y_1^n+y_2^n}\norm{e^n+e^{n-1}}$.
Taking the real part of (\ref{enplus test}), we have
$\norm{e^n+e^{n-1}}^2\leq C\varepsilon\norm{y_1^n+y_2^n}\norm{e^n+e^{n-1}}$ since $\norm{e^n+e^{n-1}}\leq C\norm{e^n+e^{n-1}}_a$.
Combining the above two inequalities, we have
\begin{equation}
\label{en1-en-1 L2 error}
\norm{e^n}^2-\norm{e^{n-1}}^2\leq 2\varepsilon\Delta t\norm{y_1^n+y_2^n}^2.
\end{equation}
Summing both sides of (\ref{en1-en-1 L2 error}) from $n=1$ to $n=N$, we have
\begin{equation}
\label{eN L2}
\norm{e^N}^2\leq 2\varepsilon\Delta t\sum_{n=1}^N\norm{y_1^n+y_2^n}^2.
\end{equation}
For $\norm{y_1^i}$, $i=1,2,...,N$, we have
\begin{align*}
\norm{y_1^i}&=\norm{\frac{\sigma(u(t_i))-\sigma(u(t_{i-1}))}{\Delta t}-\frac{u(t_i)-u(t_{i-1})}{\Delta t}}=\frac{1}{\Delta t}\norm{\int_{t_{i-1}}^{t_i}\partial_t(\sigma(u(s))-u(s))ds}\\
&\leq \frac{1}{\Delta t}\int_{t_{i-1}}^{t_i}\norm{\partial_t(\sigma(u(s))-u(s))}ds\leq \frac{1}{\Delta t}\int_{t_{i-1}}^{t_i}C\frac{H}{\Lambda\varepsilon}\cdot\varepsilon\norm{\partial^2_tu(s)}_{s^{-1}}ds\\
&\leq\frac{CH\varepsilon}{\Lambda\varepsilon}\cdot\frac{CH}{\varepsilon\min\{\varepsilon^2,\delta^2\}}=C\frac{H^2}{\Lambda\varepsilon\min\{\varepsilon^2,\delta^2\}},
\end{align*}
where Newton-Leibniz formula is used in the first line, the proof of (\ref{s norm for u-u_ms}) in Theorem \ref{convergence of the multiscale solution} has been used in the second line and  the temporal regularity in Lemma \ref{temporal regularity} has been used in the last line. Utilizing the same techniques of the above estimate, for $\norm{y_2^i}$, $i=1,2,...,N$, we have
\begin{align*}
\norm{y_2^i}&=\frac{1}{\Delta t}(\int_{t_{i-1}}^{t_{i-\frac{1}{2}}}\int_{t_{i-1}}^s\int_{t_{i-1}}^w\partial^3_tu(r)drdwds+\int_{t_{i-\frac{1}{2}}}^{t_i}\int^{t_i}_s\int^{t_i}_w\partial^3_tu(r)drdwds)\\
&\quad+\frac{1}{4}\Delta t\cdot\partial^2_tu(t_{i-1})-\frac{1}{4}\Delta t\cdot\partial^2_tu(t_i)\\
&\leq \frac{1}{\Delta t}\int_{t_{i-1}}^{t_{i-\frac{1}{2}}}(\max_{0\leq t\leq T}\norm{\partial^3_t u(t)})\cdot(s-t_{i-1})^2ds\\
&\quad+\frac{1}{\Delta t}\int_{t_{i-\frac{1}{2}}}^{t_i}(\max_{0\leq t\leq T}\norm{\partial^3_t u(t)})\cdot(s-t_{i-1})^2ds+\frac{1}{4}\Delta t\int_{t_{i-1}}^{t_i}\partial^3_t u(s)ds\\
&\leq C \Delta t^2\max_{0\leq t\leq T}\norm{\partial^3_t u(t)}\leq C\varepsilon\Delta t^2/\min\{\varepsilon^4,\delta^4\}.
\end{align*}
Then combining (\ref{eN L2}), we can estimate $\norm{e^N}$ as follows:
\begin{equation}
\label{L2 error for eN}
\norm{e^N}\leq C\sqrt{\varepsilon}(\frac{H^2}{\Lambda\varepsilon\min\{\varepsilon^2,\delta^2\}}+\frac{\varepsilon\Delta t^2}{\min\{\varepsilon^4,\delta^4\}}).
\end{equation}
Thus by (\ref{regularity for projection operator}), we obtain that
\[
\norm{u^N-u(T)}\leq \norm{e^N}+\norm{\sigma{u(T)}-u(T)}\leq C(\frac{H^2}{\Lambda\sqrt{\varepsilon}\min\{\varepsilon^2,\delta^2\}}+\frac{\varepsilon^{3/2}\Delta t^2}{\min\{\varepsilon^4,\delta^4\}}+\frac{H^2}{\Lambda\varepsilon^2}).
\]
Thus (\ref{L2 error}) is proved. Let $v=e^n-e^{n-1}$ in (\ref{en form}), then we have
\begin{equation}
\label{for energy estimate}
\frac{1}{2}(\norm{e^n}_a^2-\norm{e^{n-1}}_a^2)=i\varepsilon\frac{1}{\Delta t}\norm{e^n-e^{n-1}}^2+i\varepsilon(y_1^n+y_2^n,e^n-e^{n-1}).
\end{equation}
Taking the real part of (\ref{for energy estimate}), we have
$
\norm{e^i}_a^2-\norm{e^{i-1}}_a^2\leq 2\varepsilon\norm{y_1^i+y_2^i}\norm{e^i-e^{i-1}}
$;
Taking the imaginary part of (\ref{for energy estimate}), we have
$\norm{e^i-e^{i-1}}^2\leq \Delta t\norm{y_1^i+y_2^i}\norm{e^i-e^{i-1}}$ for $i=1,2,...,N$.
Combining the above two inequalities and
summing $i$ from 1 to $n$, we have $\norm{e^n}_a^2\leq 2\varepsilon\Delta t\sum_{i=1}^n\norm{y_1^i+y_2^i}^2$, where $n=1,2,..,N$. Similar to the proof of (\ref{L2 error for eN}), we have
\[
\norm{e^n}_a\leq C\sqrt{\varepsilon}(\frac{H^2}{\Lambda\varepsilon\min\{\varepsilon^2,\delta^2\}}+\frac{\varepsilon\Delta t^2}{\min\{\varepsilon^4,\delta^4\}}).
\]
Since $e(t)=u^n-\sigma(u(t_n))+\sigma(u(t_n))-u(t_n)$ for $t\in(t_{n-1},t_n]$, we get
\[\int_0^T\norm{e(t)}_a\leq C\sqrt{\varepsilon}(\frac{TH^2}{\Lambda\varepsilon\min\{\varepsilon^2,\delta^2\}}+\frac{\varepsilon T\Delta t^2}{\min\{\varepsilon^4,\delta^4\}}+\frac{TH}{\sqrt{\Lambda}\varepsilon}).\]

\end{proof}

\begin{remark}
If $\partial_tu(t),\partial_t^2u(t),\partial_t^3u(t)\in H^s(\Omega)$ (i.e. $\partial_tu(t),\partial_t^2u(t),\partial_t^3u(t)$ have better spatial regularity) for any $t\in[0,T]$ and $s=1,2$, we can obtain high order convergence similar to \cite{wu2022}. Here we just omit the details of the analysis under better spatial regularity assumptions.
\end{remark}
\section{Numerical experiments}\label{sec6}   
We will compare the relative errors between the numerical solution $u_{\text{cem}}$ approximated by the CEM-GMsFEM method and the reference solution $u_{\text{ref}}$ calculated by the FEM on $K_h$ or obtained by time-splitting fast spectral method \cite{BJM2002} in $L^2$ norm and $H^1$ norm with
$$\text{err}_{L^2}=\frac{\norm{u_{\text{cem}}-u_{\text{ref}}}_{L^2}}{\norm{u_{\text{ref}}}_{L^2}},\quad\text{err}_{H^1}=\frac{\norm{u_{\text{cem}}-u_{\text{ref}}}_{H^1}}{\norm{u_{\text{ref}}}_{H^1}}. $$
For simplicity, we take $\tilde{\mu}$ as
\[
  \tilde{\mu}|_{K_i}=\mu_\text{msh}\text{diam}(K_i)^{-2}|_{K_i}=12\varepsilon^2H^{-2}|_{K_i}
\]
for all numerical experiments, as suggested in \cite{Ye2023}.
%In the following discussions, we mark statements of direct observations from numerical experiments with a circled number, e.g., \circled{1}.
We implement the method using the Python libraries NumPy and SciPy.
Moreover, we will show the performance of our method for the approximation of observables, including the position density
\begin{equation}\label{pos}
n(x,t)=|u(x,t)|^2,    
\end{equation}
and the energy density
\begin{equation}\label{ene}
e(x,t)=\frac{\varepsilon^2}{2}|\nabla u(x,t)|^2+V(x)|u(x,t)|^2.    
\end{equation}
\subsection{1D potential}
\begin{figure}[!ht]
\centering
\begin{subfigure}{0.45\textwidth}
\includegraphics[width=\linewidth]{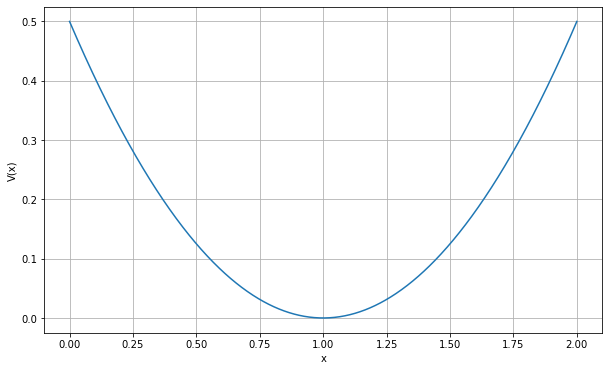}
\caption{$1D$ smooth Potential}
\label{fig:1}
\end{subfigure}
\begin{subfigure}{0.45\textwidth}
\includegraphics[width=\linewidth]{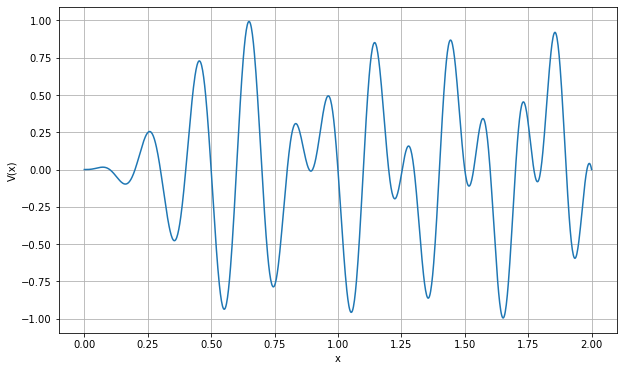}
\caption{ 1D multiplicative two-scale potential }
\label{fig:2}
\end{subfigure}
\caption{1D Potential}
\label{fig:3}
\end{figure}

In this section, we present our methods in 1D smooth potential and  1D multiplicative two-scale potential as shown in \cref{fig:3}. We conside the final time $T=0.1$ and domain $\Omega=[0,2].$ The $1D$ initial data is chosen as 
$$u_{\text{int}} =\sqrt{r_0}\text{exp}(\mathrm{i}S_0 / \varepsilon),$$
where $$r_0 = (\text{exp}(-50 (x - 1) ^2))^ 2, S_0 = -0.2\text{log}(\text{exp}(5(x - 1)) + \text{exp}(-5(x - 1))). $$
\subsubsection{1D smooth potential}
\begin{figure}[!ht]
\centering
\includegraphics[width=\linewidth]{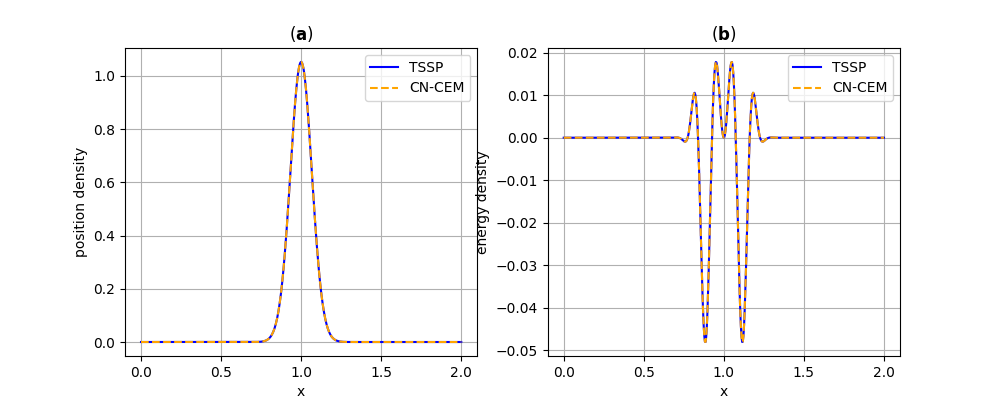}
\caption{when $\varepsilon=1/32$. Left: Profiles of numerical and reference position density functions;  Right: Profiles of numerical and reference energy density functions. }
\label{fig:4}
\end{figure}

\begin{figure}[!ht]
\centering
\includegraphics[width=\linewidth]{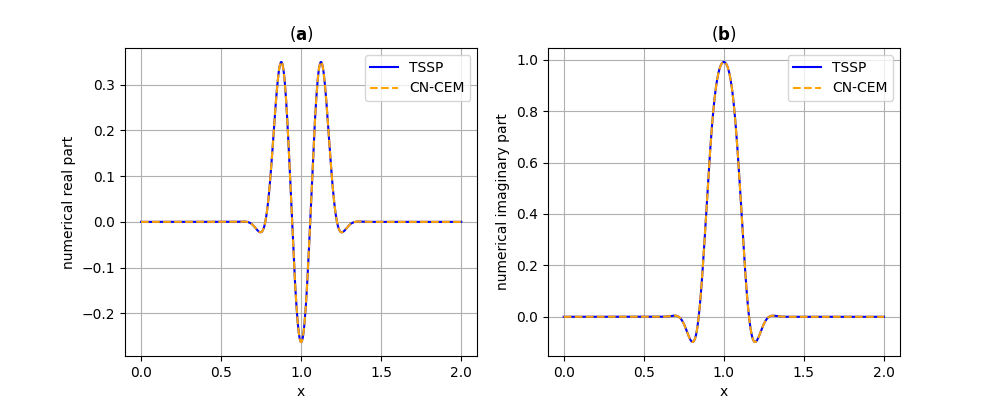}
\caption{Real and imaginary parts of the wavefunction $u(t,x)$ when $t=1$ and $\varepsilon=1/32$. Left: real part of the numerical solution;  Right: imaginary part of the numerical solution. }
\label{fig:5}
\end{figure}

The potential is defined in the following and shown in \cref{fig:1}
$$V(x)=0.5(x-1)^2.$$
\circled{1} {$\varepsilon=1/32, H=1/64, \Delta t=10^{-2}$.} Here, we use the Time splitting spectral method (TSSP) solution  as the reference solution, and compare the relative errors between the Crank-Nicolson CEM scheme and reference solution. The oversampling size $m$ of the multiscale scheme is chosen as $m=\frac{2}{3}\log_2(2/H)$, and the numerical results of two density functions defined in \cref{pos,ene} are depicted in \cref{fig:4}, showing a good approximation of our multiscale method. In \cref{fig:5}, we show the real part and imaginary part of the wave function to show the oscillation phenomenon when $\varepsilon=1/32$. \circled{2} According to the theoretical result of Theorem \ref{final convergence}, the error is expected to increase with the decrease of $\varepsilon$. However, despite this theoretical prediction, the observed $L^2$ and $H^1$ errors remain stable in the order of $O(10^{-4})$  in Table \ref{1dtab}. This indicates that our new methods demonstrate significant robustness of managing error growth even as $\varepsilon$ decreases. {The reason is that we choose $ H/\sqrt{\Lambda}=O(\varepsilon^{\frac{5}{4}}), \Delta t=O(\varepsilon^{\frac{5}{4}})$ based on the result of \cref{final convergence}.}
%and the initial function is 
%$$u_{\text{int}} =\sqrt{r_0}\text{exp}(\mathrm{i}S_0 / \varepsilon).$$
%where
%$$r_0 = (\text{exp}(-50 (x - 1) ^2))^ 2, S_0 = -0.2\text{log}(\text{exp}(5(x - 1)) + \text{exp}(-5(x - 1))). $$
%We choose the following parameters:
%$$dt=0.01,\quad \text{timestep}=100,\quad L=2,\quad N=2^8. $$
%We use the Time splitting spectral method (TSSP) solution as the reference solution, and compare the differences between Crank-Nicolson CEM schemes with respect to the reference solution. The {oversampling size} for the multiscale scheme is chosen as $l=\log_2(1/H)$ which can be obtained from the Assumption 3 in theoretical analysis. We choose $H$ to be $1/(2^3),$ and the numerical results of two density functions defined in \cref{pos,ene} are depicted in \cref{fig:1d1}. In \cref{fig:1D3}, we show the real part and imaginary part of the wave function to show the oscillation phenomenon when $\varepsilon=1/32$. According to the theoretical result of Theorem \ref{final convergence}, the error is expected to increase with the decrease of $\varepsilon$. However, despite this theoretical prediction, the observed $L^2$ and $H^1$ errors remain stable in the order of $O(10^{-4})$  in Table \ref{1dtab}. This indicates that our new methods demonstrate significant robustness, effectively managing error growth even as $\varepsilon$ decreases.
\begin{table}[htbp]
\footnotesize
\caption{Relative $L_2$ and $H^1$ Errors with different $\varepsilon$.}\label{1dtab}
\begin{center}
  \begin{tabular}{|c|c|c|} \hline
$\varepsilon$ & $L^2$ Error &$H^1$ Error \\ \hline
    $1/128$ & $3.338e-4$ & $7.828e-4$ \\
$1/64$   & $1.754e-4$ & $4.076e-4$ \\ 
$1/32$& $1.558e-4$ &$4.015e-4$  \\ \hline
  \end{tabular}
\end{center}
\end{table}
\subsubsection{1D multiplicative two-scale potential}
\label{1D multiplicative two-scale potential example}
\begin{figure}[!ht]
\centering
\includegraphics[width=\linewidth]{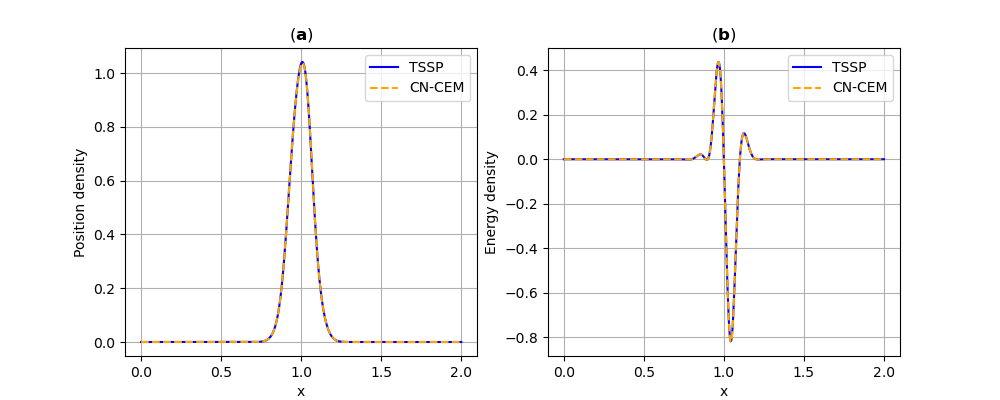}
 \caption{$\varepsilon=1/16$. Left: Profiles of numerical and reference position density functions;  Right: Profiles of numerical and reference energy density functions.}
    \label{fig:1d-hc1}
\end{figure}

\begin{figure}[!ht]
\centering
\includegraphics[width=\linewidth]{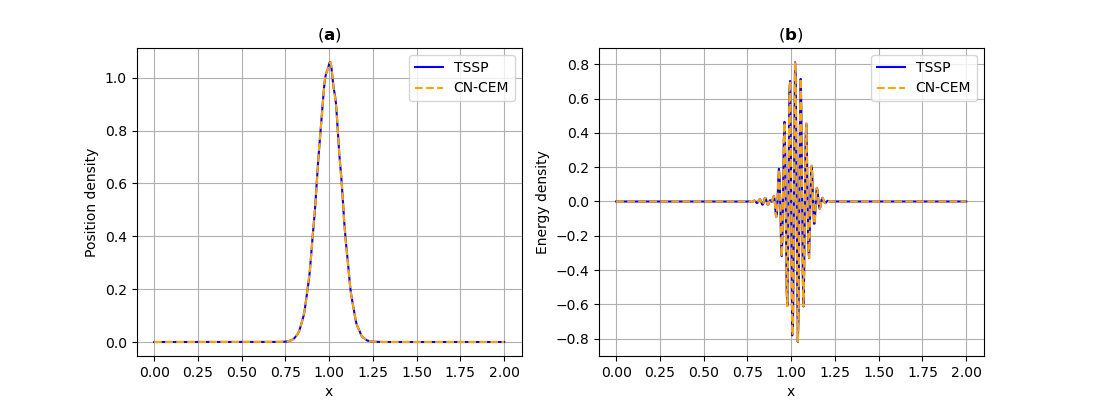}
\caption{ $\varepsilon=1/32$. Left: Profiles of numerical and reference position density functions; Right: Profiles of numerical and reference energy density functions.}
    \label{fig:1de2}
\end{figure}

In this example, we choose the potential to be the following and plot it in  \cref{fig:2}. 
$$V(x)=\sin(\frac{x^2}{\delta_1})\sin(\frac{\pi x}{\delta_2}),\quad x\in [0,2].$$
\circled{1} {$(\varepsilon\leq \delta) \varepsilon=1/16,\delta_1=1/4, \delta_2=1/10$,} and we use the Time splitting spectral method (TSSP) solution as $u_{\text{ref}}$, and the solution of  the Crank-Nicolson CEM schemes to be $u_{\text{cem}}$. The oversampling size $m$ of the multiscale scheme is chosen as $m=\frac{2}{3}\log_2(2/H)$ which can be obtained from the Assumption 3 in theoretical analysis. We choose $H$ to be $1/64,$ and {the position density and energy density
of the reference solution and the multiscale solution  are depicted in \cref{fig:1d-hc1}, showing a good approximation.} In \cref{tab:1d-hc3}, we also observe $L^2$ and $H^1$ errors remain stable in the order of $O(10^{-1})$. This indicates that our new methods demonstrate significant robustness to  maintain the error growth even as $\varepsilon$ decreases. \circled{2} {$(\varepsilon\geq\delta) \varepsilon=1/32,\delta_1=1/4, \delta_2=1/64$, the numerical results of two density functions defined in \cref{pos,ene} are depicted in \cref{fig:1de2}, also demonstrating a good approximation. However, for the smaller \(\varepsilon\), there is a larger magnitude of the energy density functions and more oscillations compared with $\varepsilon=1/16$ in \cref{fig:1d-hc1}.}
\begin{table}[!ht]
\footnotesize
\caption{Relative $L_2$ and $H^1$ Errors with different $\varepsilon$.}
\begin{center}
\begin{tabular}{|c|c|c|} \hline
$\varepsilon$ & $L^2$ Error &$H^1$ Error \\ \hline
$1/128$ & $6.548e-01$ & $8.131e-01$ \\
$1/64$   & $2.336e-01$ & $4.095e-01$ \\ 
$1/32$& $1.138e-01$ &$2.022e-01$  \\ \hline
\end{tabular}
\end{center}
\label{tab:1d-hc3}
\end{table}
\subsection{2D potential}
\begin{figure}[htbp]
\centering
\begin{subfigure}{0.45\textwidth}
\includegraphics[width=\linewidth]{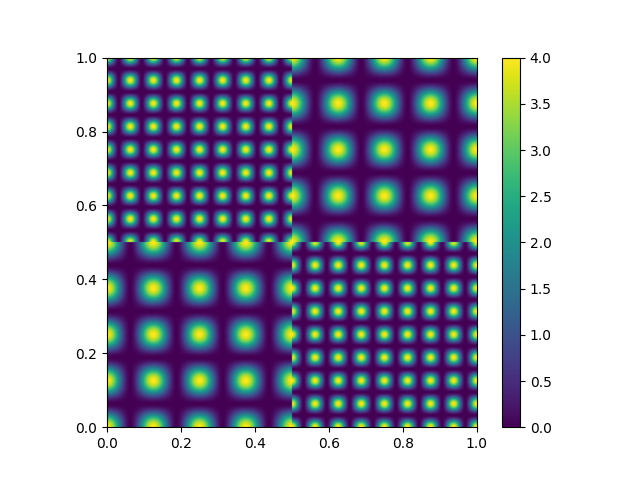}
\caption{$2D$ checkboard Potential}
\label{fig:medium1}
\end{subfigure}
\begin{subfigure}{0.45\textwidth}
\includegraphics[width=\linewidth]{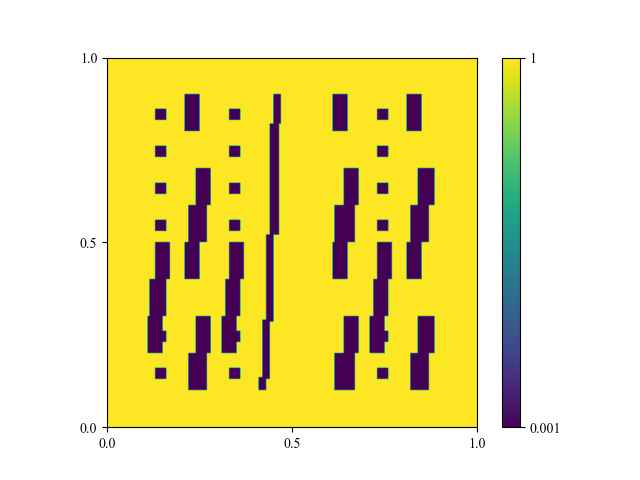}
\caption{$2D$ high-contrast Potential}
\label{fig:source1}
\end{subfigure}
\caption{2D Potential}
\label{2D}
\end{figure}

In this section, we present our methods in 2D checkboard potential and heterogeneous potential with high contrast, as shown in \cref{2D}. We consider the final time $T=1$ and domain $\Omega=[0,1]^2.$ The $2D$ initial data is chosen as 
{
$$u_0(x,y)=(\frac{10}{\pi})^{1/2}\text{exp}(-5(x-1/2)^2-5(y-1/2)^2)\text{exp}(-\mathrm{i}\frac{(x-1/2)^2+(y-1/2)^2)}{\varepsilon}).$$}
\subsubsection{2D checkboard potential}
\begin{figure}[!ht]
\centering
\includegraphics[width=\linewidth]{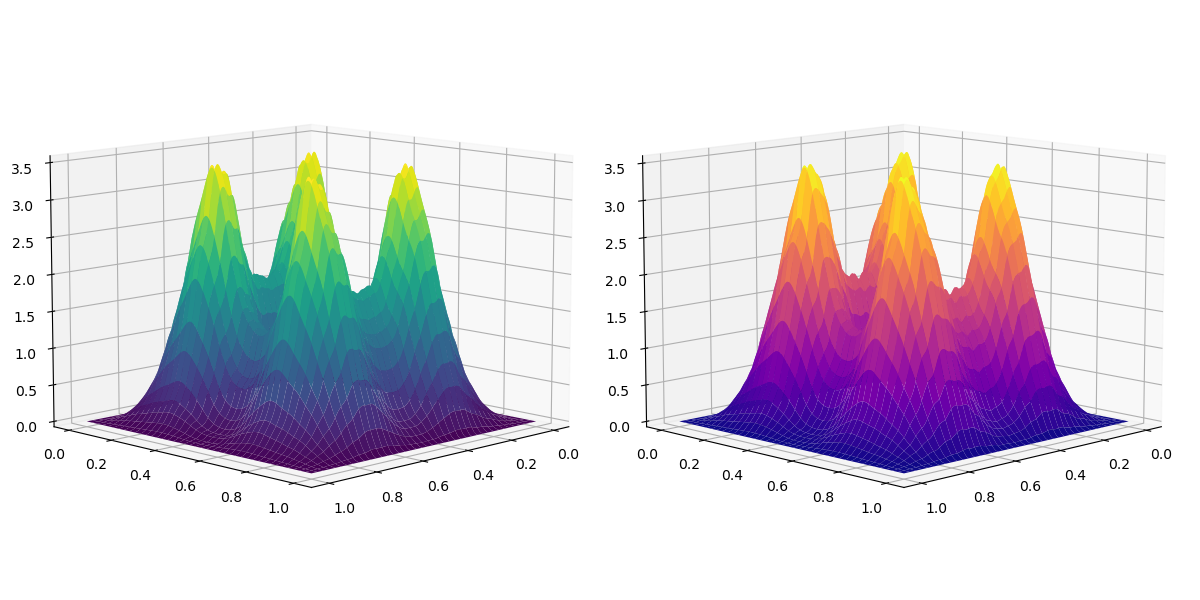}
\caption{ $\varepsilon=1/8$ and $\delta=1/16$ . Left: Profiles of CN-FEM position density functions; Right: Profiles of CN-CEM position functions.}
\label{fig:checkpot1}
\end{figure}

\begin{figure}[!ht]
\centering
\includegraphics[width=\linewidth]{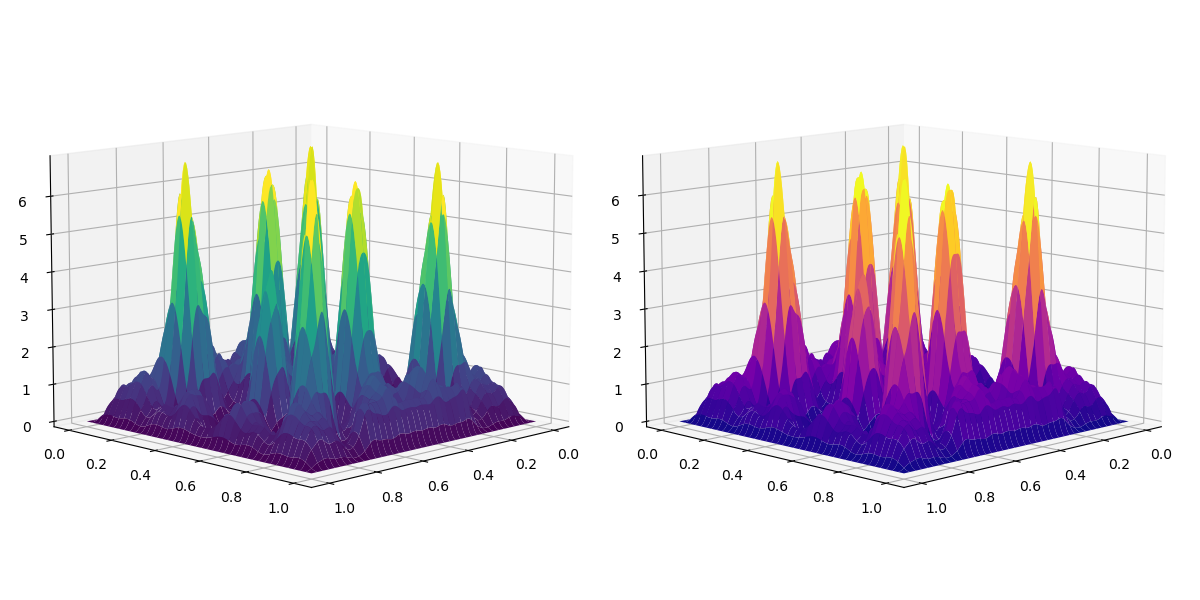}
\caption{ $\varepsilon=1/32$ and $\delta=1/16$.  Left: Profiles of CN-FEM  position density functions; Right: Profiles of CN-CEM position density functions.}
\label{fig:checkpot2}
\end{figure}

The checkboard potential contains different
lattice structures and admits discontinuities along the interfaces as in quantum metamaterials
\cite{Quach11}.
\begin{equation*}\label{check}
 V=
\begin{cases}
(\cos(\frac{2\pi x}{\delta_2})+1)  (\cos (\frac{2\pi y}{\delta_2})+1), &\{0\leq x,y\leq 1/2 )\cup\{1/2 \leq x,y\leq 1\},\\
(\cos(\frac{2\pi x}{\delta_1})+1)(\cos (\frac{2\pi y}{\delta_1})+1),&\text{otherwise}.
    \end{cases}
\end{equation*}
\circled{1} {$(\varepsilon\geq\delta)\, \varepsilon=1/8, \delta_1=1/8, \delta_2=1/16.$ } The reference solution $u_{\text{ref}}$ is computed by the Crank-Nicolson FEM scheme. We fixed the fine grid to be 200 and the coarse grid is chosen within $  1/10, 1/20, 1/40$. Due to the relationship $m=\frac{2}{3}\log_2(1/H)$ between the {oversampling size $m$} and the coarse mesh $H$ from the Assumption 3 at theoretical analysis, we choose the {oversampling size $m$} to be $=2,3,4$. We also fix the $l_j=3$, indicating that we calculate the first three eigenfunctions and construct three multiscale bases for each coarse element. {In \cref{tab:2d1}, we achieve the  suggested convergence order of both the $L^2$ norm and energy norm to be 2 and 1, demonstrating the stability of our newly constructed CN-CEM scheme with the $\varepsilon$-dependent initial data.} \circled{2} {$(\varepsilon\leq\delta)\,\varepsilon=1/32, \delta_1=1/8, \delta_2=1/16$.} The fine mesh is chosen to be $1/200$, and the coarse mesh is chosen within  $  1/10, 1/20, 1/40$. In \cref{tab:2d5}, we have  convergence rates 2 and 1 in $L^2$ norm and $H^1$ norm, respectively. \circled{3} Additionally, when we choose $\Delta t=1/2^5, H=1/40, \varepsilon=1/16$, {\cref{tab:2d2} presents the relative errors of CN-CEM schemes in both the $L^2$ and $H^1$ norms, indicating that the errors keep stable as the value of $\delta$ decreases, aligning with our theoretical results in \cref{final convergence}.} \circled{4} {We visualize the position density functions of two schemes in \cref{fig:checkpot1} with $\varepsilon=1/8$ and \cref{fig:checkpot2} with $\varepsilon=1/32$, which reveal higher oscillations of smaller $\varepsilon$  and good approximations of multiscale methods.}
\begin{table}[!ht]
\centering
\caption{Relative errors in the $L^2$ norm and energy norm with different $H$.}
\begin{tabular}{ c c c c c c c }
\hline\hline
\multicolumn{5}{c}{$(\varepsilon\geq\delta) \quad \varepsilon=1/8,\,\delta=1/16,\, t=1,\, l_j=3$} \\
\hline
$H$ & $h$ & $m$ & $\|u_{\text{cem}}-u_{h}\|_{L^2(\Omega)}$ & Order & $\|u_{\text{cem}}-u_{h}\|_{H^1(\Omega)}$ & Order \\[0.5ex]
\hline
1/10 & 1/200 & 2 & 3.856e-2 &  &1.651e-1  &  \\
1/20 & 1/200  & 3 & 8.717e-3 & 2.14 & 7.067e-2 & 1.22\\
1/40 & 1/200  & 4 & 1.831e-3 &2.25 &2.325e-2 &1.60 \\
\hline
\end{tabular}
\label{tab:2d1}
\end{table}
\begin{table}[!ht]
\centering
\caption{Relative errors in the $L^2$ norm and energy norm with different $H$.}
\begin{tabular}{ c c c c c c c }
\hline\hline
\multicolumn{5}{c}{$(\varepsilon\leq\delta)\quad \varepsilon=1/32,\, \delta=1/16,\, t=1,\, l_j=3$} \\
\hline
$H$ & $h$ & $m$ & $\|u_{\text{cem}}-u_{h}\|_{L^2(\Omega)}$ & Order & $\|u_{\text{cem}}-u_{h}\|_{H^1(\Omega)}$ & Order \\[0.5ex]
\hline
1/10 & 1/200 & 2 & 3.397e-2 &  &  6.955e-2&  \\
1/20 & 1/200  & 3 & 7.413e-3& 2.19 & 3.963e-2& 0.81\\
1/40 & 1/200  & 4 &  1.803e-3& 2.03 & 1.423e-2&1.48 \\
\hline
\end{tabular}
\label{tab:2d5}
\end{table}
\begin{table}[!ht]
\centering
\caption{Relative errors in the $L^2$ norm and energy norm with different $\varepsilon$.}
\begin{tabular}{ c c c c c }
\hline\hline
\multicolumn{5}{c}{$\varepsilon=1/16, H=1/40, t=1, l_j=3$} \\
\hline
$\delta$ & $h$ & $m$ & $\|u_{\text{cem}}-u_{h}\|_{L^2(\Omega)}$ & $\|u_{\text{cem}}-u_{h}\|_{H^1(\Omega)}$ \\[0.5ex]
\hline
1/8 & 1/200 & 4 & 2.434e-2  & 2.267e-1 \\
1/16& 1/200  & 4 & 4.768e-2 & 2.767e-1\\
1/32& 1/200  & 4 & 7.589e-2  & 3.386e-1 \\
\hline
\end{tabular}
\label{tab:2d2}
\end{table}
\subsubsection{2D High-contrast potential}
\label{2D High-contrast potential example}
\begin{figure}[!ht]
\centering
\includegraphics[width=\linewidth]{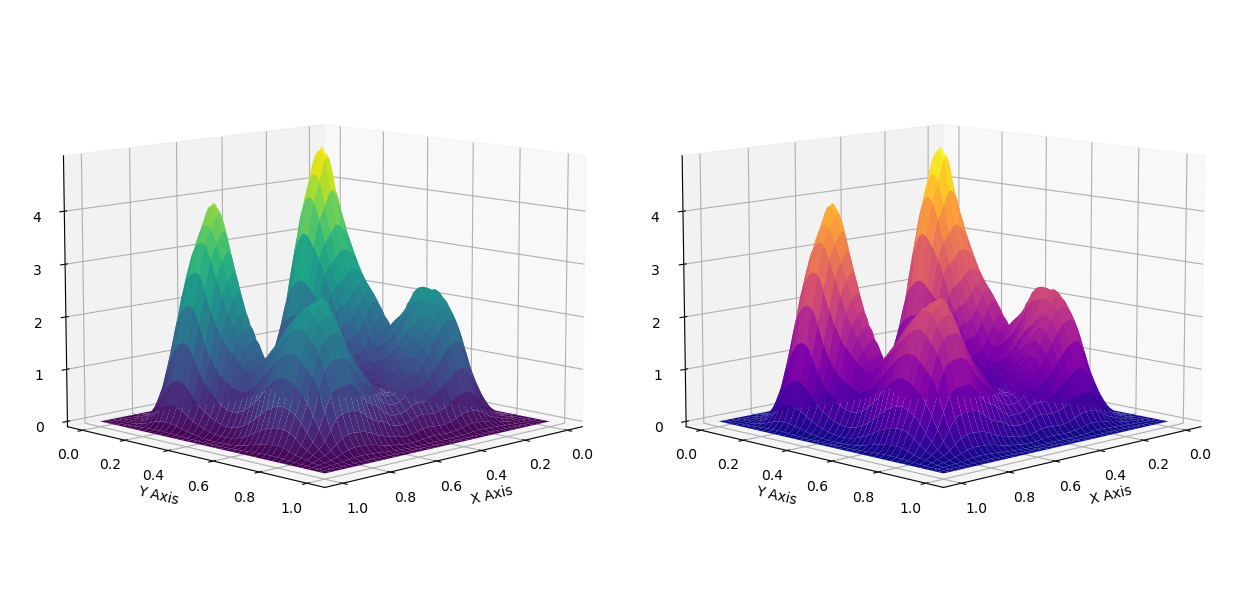}
\caption{ $\Upsilon=10^3$ and $\varepsilon=1/8$. Left: Position density of CN-FEM solution; Right: Position density of the CN-CEM solution.}
\label{fig:2dtop}
\end{figure}

In this subsection, we consider a random inclusion model that is utilized in several multiscale methods as a showcase of the ability to handle non-periodic coefficient profiles \cite{jin2024,ye2024, Chung2018}. The high contrast ratio is chosen as $\Upsilon=10^{3}$ where $V(x,y)$ is taken 1 or $10^{-3}$ in the domain and displayed in \cref{2D} (b). The reference solution $u_\text{ref}$ is  computed by the Crank-Nicolson FEM scheme. 
\begin{table}[!ht]
\centering
\caption{Relative errors in the $L^2$ norm and energy norm with different $H$.}
\begin{tabular}{ c c c c c }
\hline\hline
\multicolumn{5}{c}{$\Upsilon=10^3, \varepsilon=1/8, t=1,  l_j=3$} \\
\hline
$H$ & $h$ & $m$ & $\|u_{\text{cem}}-u_{h}\|_{L^2(\Omega)}$ & $\|u_{\text{cem}}-u_{h}\|_{H^1(\Omega)}$ \\[0.5ex]
\hline
1/10 & 1/200 & 2 & 3.343e-3  & 2.309e-2 \\
1/20 & 1/200  & 3 & 5.720e-4 & 8.445e-3\\
1/40 & 1/200  & 4 & 3.374e-4  & 7.569e-3 \\
\hline
\end{tabular}
\label{tab:2d3}
\end{table}
\begin{figure}[!ht]
\centering
\includegraphics[width=\linewidth]{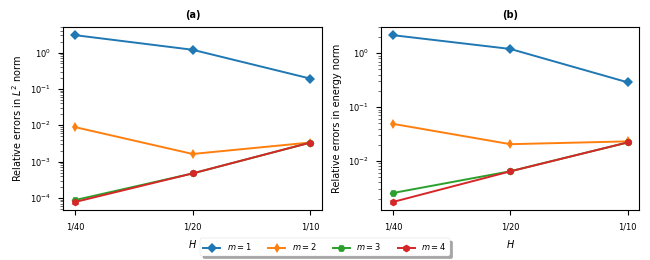}
\caption{$\Upsilon=10^3$ and $\varepsilon=1/8$. (a) Relative Error in $L_2$ . (b) Relative $H_1$ Norm with different {oversampling size} and coarse mesh $H$ when $t=1$.}
\label{fig:2derror}
\end{figure}
\circled{1}  We fixed $\varepsilon=1/8,\Delta t=1/2^5$, the fine grid to be $1/200$ and the coarse grid is chosen within $1/10, 1/20, 1/40$, and we choose the {oversampling size $m$} to be $1,2,3,4$. Based on Assumption 3, it is clear that as the high contrast increases, a smaller oversampling is preferable. We also fix the $l_j$=3, indicating that we calculate the first three eigenfunctions and construct three multiscale bases for each coarse element. { Although we our initial data is $\varepsilon$-dependent, we plot and show the relative error in \cref{fig:2derror} and \cref{tab:2d3} by adjusting the oversampling layers $m=\{1,2,3,4\}$, and it is obvious that when the {oversampling size} equals to 4 we can have the $H$ convergence of the energy norm in \cref{fig:2derror} (b). In this example, we can see that the proposed CN-CEM numerical method can capture the high-contrast potential and the small $\varepsilon$ robustly , which also shows that the choosing of suitable $m$ is very important.} \circled{2} {In order to show our CN-CEM scheme is also robust for the high-contrast case, we vary the high contrast ratios $\Upsilon=\{10,10^2,10^3,10^4\}$ of the potential and give the relative  $L^2$ and $H^1$ errors in  \cref{tab:2d4}. It is encouraging to see that the relative $L^2$ error remains in the order of $O(10^{-4})$ and the relative $H^1$ error remains in the order of $O(10^{-3})$ when we fixed $\varepsilon$, the coarse mesh $H$, {oversampling size}.} \circled{3} In \cref{fig:2dtop}, We visualize the position density of the CN-FEM and CN-CEM schemes, highlighting the oscillations that occur in regions of high contrast $\Upsilon=10^3$, which also shows good approximation of our new scheme.
\begin{table}[!ht]
\footnotesize
\caption{Relative $L^2$ and $H^1$ Errors with different contrast ratios when $\varepsilon=1/8$.} 
\begin{center}
\begin{tabular}{|c|c|c|c|c|} \hline
$\Upsilon$ & $H$ & $m$ & $L^2$ Error & $H^1$ Error \\ \hline
$10$ & 1/40& 4 &  1.817e-4  &     5.605e-3  \\
$10^2$ & 1/40& 4  &  2.318e-4   &  6.480e-3  \\ 
$10^3$ & 1/40& 4  &  3.374e-4    &   7.569e-3     \\ 
$10^4$ & 1/40&  4  &  4.380e-4      & 8.577e-3      \\ \hline
  \end{tabular}
\end{center}
\label{tab:2d4}
\end{table}
\section{Conclusions}\label{sec7}
In this paper, we introduce a novel method for solving the Schr\"{o}dinger equation with high-contrast potentials in the semiclassical regime, inspired by the newly proposed multiscale method known as CEM-GMsFEM. This multiscale approach leverages the LOD method by constructing new multiscale basis functions to approximate solutions through oversampling techniques in a more generalized manner. In the analysis section, we first establish the regularity of the spatial solution. We present a significant result showing that the error estimate between the multiscale solution and the reference solution is independent of the parameters $\varepsilon,\delta, H,\Lambda$, thanks to Assumption \ref{oversampling layers} concerning the {oversampling size}. Next, we address the temporal regularity of our newly constructed CN-CEM scheme. The final theorem, detailed in Theorem \ref{final convergence}, indicates that the right-hand side of the error estimate remains bounded by a small number when considering the relationship between the coarse mesh size $H$, $\Delta t$ and $\varepsilon,\delta$. In the numerical section, we demonstrate effective approximations for the 1D potential, illustrating both the real and imaginary parts of the CN-CEM solution. For the 2D checkerboard example, we achieve second-order accuracy in the relative $L^2$ norm and first-order accuracy in the relative $H^1$ norm, consistent with our theoretical predictions across various {oversampling sizes} and coarse mesh sizes $H$. Finally, for high-contrast potentials, we explore scenarios with both large and small 
$\varepsilon$. Our results confirm the anticipated relationships among the high-contrast ratio, coarse mesh size, and $\varepsilon$. This is the first investigation into high-contrast potentials in the context of the Schr\"{o}dinger equations.

In the future, we plan to extend the CEM-GMsFEM to address more general types of potentials, including time-dependent and random potentials. Additionally, we are exploring the integration of this multiscale method with deep learning frameworks to tackle a broader range of problems in kinetic fields and uncertainty quantification.

\bibliographystyle{siamplain}
\bibliography{References}
\end{document}